\documentclass[11pt]{amsart}
\pdfoutput=1
 \usepackage[colorlinks=true, pdfstartview=FitV, linkcolor=blue, citecolor=blue, urlcolor=blue]{hyperref}
 
 \usepackage{graphicx}
 
\newcommand{\aaa}{{\mathcal A}}

\newcommand{\lll}{{\mathcal L}}
\newcommand{\ppp}{{\mathcal P}}
\newcommand{\sss}{{\mathcal S}}

\newcommand{\T}{{\mathcal T}}
\newcommand{\R}{{\mathbb R}}

\newcommand{\Z}{{\mathbb Z}}
\newcommand{\C}{{\mathbb C}}

\newcommand{\vecv}{{\vec{v}}}

\newtheorem{thm}{Theorem}[section]
\newtheorem*{thm*}{Theorem}

\newtheorem{dfn}[thm]{Definition}

\everymath{\displaystyle}

\theoremstyle{definition}
\newtheorem{ex}{Example}

\begin{document}
 \title{A primer on substitution tilings of the Euclidean plane}
 \author{Natalie Priebe Frank }

\address{Natalie Priebe Frank\\Department of Mathematics\\Vassar College\\Box 248\\Poughkeepsie, NY  12604\\nafrank@vassar.edu}

\maketitle

\begin{abstract}
This paper is intended to provide an introduction to the theory of substitution tilings.   
For our purposes, tiling substitution rules are divided into two broad classes: geometric and combinatorial.   Geometric substitution tilings include self-similar tilings such as the well-known Penrose tilings;  for this class there is a substantial body of research in the literature.   Combinatorial substitutions are just beginning to be examined, and some of what we present here is new.
We give numerous examples, mention selected major results, discuss connections between the two classes of substitutions, include current research perspectives and questions, and provide an extensive bibliography.
Although the author attempts to fairly represent the as a whole, the paper is not an exhaustive survey, and she apologizes for any important omissions.
\end{abstract}


\section{Introduction}
A {\em tiling substitution rule} is a rule that can be used to construct infinite tilings of $\R^d$ using a finite number of tile types.  The rule tells us how to ``substitute" each tile type by a finite configuration of tiles in a way that can be repeated, growing ever larger pieces of tiling at each stage.    In the limit, an infinite tiling of $\R^d$ is obtained.

In this paper we take the perspective that there are two major classes of tiling substitution rules:  those based on a linear expansion map and those relying instead upon a sort of ``concatenation" of tiles.   The first class, which we call {\em geometric tiling substitutions}, includes {\em self-similar tilings}, of which there are several well-known examples including the Penrose tilings.   In this class a tile is substituted by a configuration of tiles that is a linear expansion of itself, and this geometric rigidity has permitted quite a bit of research to be done.  We will note some of the fundamental results,  directing the reader to appropriate references for more detail.
The second class, which we call {\em combinatorial tiling substitutions}, is sufficiently new that it lacks even an agreed-upon definition.   In this class the substitution rule replaces a tile by some configuration of tiles that may not bear any geometric resemblance to the original. 
The difficulty with such a rule comes when one wishes to iterate it: we need to be sure that the substitution can be applied repeatedly so that all the tiles fit together without gaps or overlaps.   The examples we provide are much less well-known (in some cases new) and are ripe for further study.    The two classes are related in a subtle and interesting way that is not yet well understood.

\subsection{Some history}
The study of aperiodic tilings in general, and substitution tilings specifically, comes from the confluence of several discoveries and lines of research.   Interest in the subject from a philosophical viewpoint came to the forefront when Wang \cite{Wang} asked about the decidability of the ``tiling problem": whether a given set of prototiles can form an infinite tiling of the plane.   He tied this answer to the existence of  ``aperiodic prototile sets": finite sets of tiles that can tile the plane, but only nonperiodically.    He saw that the problem is undecidable if an aperiodic prototile set exists.  Berger \cite{Berger} was the first to find an aperiodic prototile set and was followed by many others, including Penrose.   It turned out that one way prove a prototile set is aperiodic involves showing that every tiling formed by the prototile set is self-similar.

Independently, work was proceeding on one-dimensional symbolic substitution systems, a combination of dynamical systems and theoretical computer science.   Symbolic dynamical systems had become of interest due to their utility in coding more complex dynamical systems, and great progress was being made in our understanding of these systems.   Queffelec \cite{Queffelec} summarized what was known about the ergodic and spectral theory of substitution systems, while a more recent survey of the state of the art appears in \cite{Fogg}.   Substitution tilings can be seen as a natural extension of this branch of dynamical systems; insight and proof techniques can often be borrowed for use in the tiling situation.   We will use symbolic substitutions motivate our study in the next section.

From the world of physics, a major breakthrough was made in 1984 by Schechtman et. al. \cite{Schechtman} with the discovery of a metal alloy that, by rights, should have crystalline since its x-ray spectrum was diffractive.   However, the diffraction pattern had five-fold rotational symmetry, which is not allowed for ideal crystals!   This type of matter has been termed ``quasicrystalline", and self-similar tilings like the Penrose tiling, having the right combination of aperiodicity and long-range order, were immediately recognized as valid mathematical models.   Dynamical systems entered the picture, and it was realized that the spectrum of a tiling dynamical system is closely related to the diffraction spectrum of the solid it models \cite{Dworkin,Hof}.  Thus we find several points of departure for the study of substitution tilings and their dynamical systems.

\subsection{One-dimensional symbolic substitutions}
\label{one.d.subs}

Let $\aaa$ be a finite set called an {\em alphabet}, whose elements are called {\em letters}.   Then $\aaa^*$, the set of all finite {\em words} with elements from $\aaa$, forms a semigroup under concatenation.   A {\em symbolic substitution} is any map $\sigma: \aaa \to \aaa^*$.   A symbolic substitution can be applied to words in $\aaa^*$ by concatenating the substitutions of the individual letters.   A block of the form $\sigma^n(a)$ will be called a {\em level-$n$ block of type $a$}.

\begin{ex}
\label{const1D}
Let $\aaa = \{a, b\}$ and let $\sigma(a)= a b b$ and $\sigma(b) = bab$.   Beginning with the letter $a$ we get
$$a \to abb \to abb\, bab \, bab \to abb \,bab\, bab \,\,bab \,abb\, bab \,\,bab\, abb\, bab \to \cdots,$$
where we've added spaces to emphasize the breaks between substituted blocks.   Notice that the block lengths
triple when substituted.
\end{ex}

\begin{ex}
\label{nonconst1D}
Again let $\aaa = \{a,b\}$; this time let $\sigma(a) = a b$ and $\sigma(b) = a$.   If we begin with $a$ we get:
$$a \to ab \to ab \, a \to ab \, a \,\, ab \to ab \, a \,\, ab \,\,ab \, a \to ab \, a \, ab \,\,ab \, a\,\,\,  ab \, a \,\, ab \to \cdots$$
Note that in this example block lengths are 1, 2, 3, 5, 8, 13, ... , and the reader can verify that they will continue growing as Fibonacci numbers.
\end{ex}

These examples illustrate the major distinction we make between substitutions.   In the first example, the length of a substituted letter is always 3 and thus the size of any level-$n$ block must be $3^n$;  this is a {\em substitution of constant length}.   In the second example the size of a substituted letter depends on the letter itself, and the size of a level-$n$ block is computed recursively; this is a {\em substitution of non-constant length}.   This is the essence of the distinction between geometric and combinatorial tiling substitutions.

It is interesting to consider infinite sequences of the form $\{x_k\} = ...x_{-2}x_{-1}.x_0x_1x_2...$  in $\aaa^{\Z}$.   Such a sequence is said to be {\em admitted} by the substitution if every finite block of letters is contained in some level-$n$ block.  In the theory of dynamical systems, the space of all sequences admitted by the substitution is studied using the shift action $s(\{x_k\} )= \{x_{k+1}\}$ (basically, moving the decimal point one unit to the right). An interested reader should see \cite{Queffelec, Fogg} to find out more.

\subsection{Two-dimensional symbolic substitutions}
The most straightforward generalization to tilings of $\R^2$ (or $\R^d$) is to use unit square tiles labeled (colored) by the alphabet $\aaa$.   These tilings can be considered as sequences in $\Z^2$, and substitutions can take letters to square or rectangular blocks of letters.   We only need to ensure that all of the blocks ``fit" to form a sequence without gaps or overlaps.

The constant length case is to expand each colored tile by some integer $n>1$ and then subdivide into $n^2$ (or $n^d$) colored unit squares.    A simple method for the non-constant length case is to take the direct product of one-dimensional substitutions of non-constant length.  

\begin{ex} \label{const2D}
Let $\aaa= \{1,2\}$, where we represent $1$ as a white unit square tile and $2$ as a blue unit square tile.  Suppose the length expansion is 3 and that the tiles are substituted by a three-by-three array of
tiles, colored  as in Figure \ref{square1}.  Starting with the blue level-0 tile, level-0, level-1, level-2 and level-3 tiles are shown in Figure \ref{square2}.   One sees in this example the tiling version of the rule creating the Sierpinski carpet.

\begin{figure}[ht]
\includegraphics[width=2in]{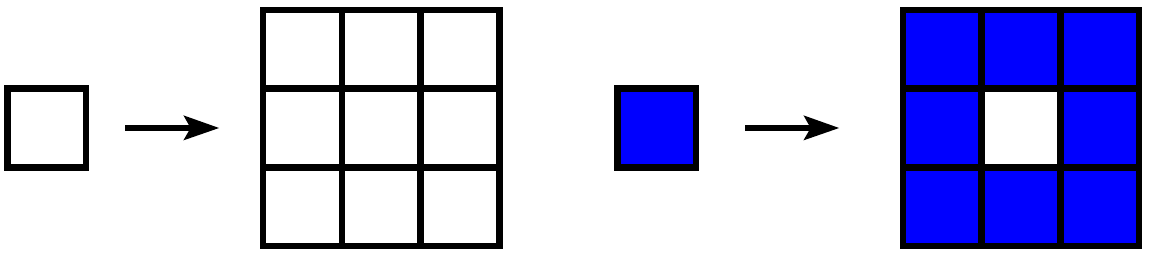}
\caption{A substitution on two colored square tiles.}
\label{square1}
\end{figure}

\begin{figure}[ht]
\includegraphics[width=5in]{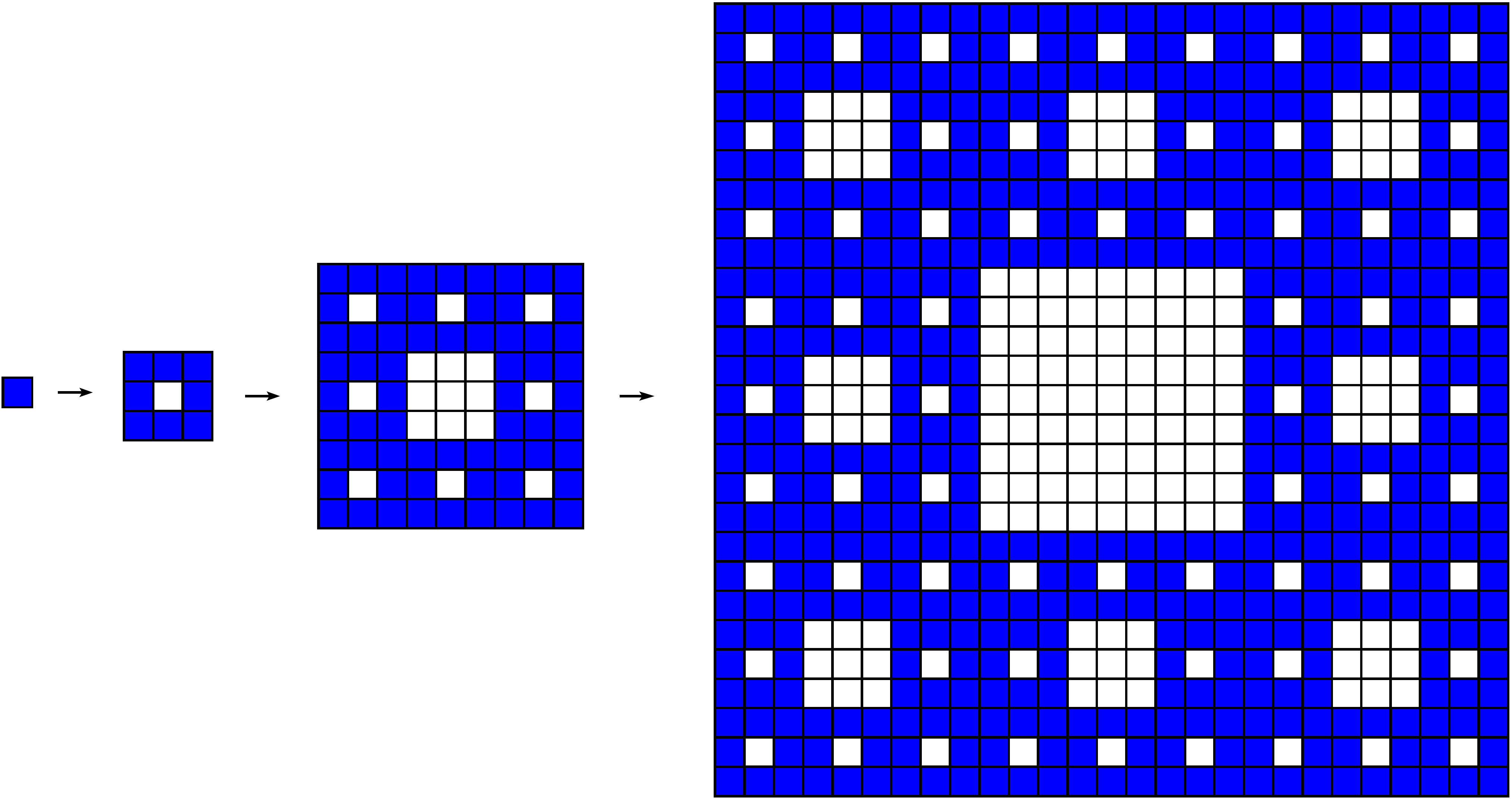}
\caption{Level-0, level-1, level-2, and level-3 tiles.}
\label{square2} 
\end{figure}
\end{ex}

\begin{ex}\label{nonconst2D}
This time, let the alphabet be $\{a,b\} \times \{a,b\} $; for simplicity of notation we put $(a,a) = 1, (a,b) = 2, (b,a) = 3, (b,b) = 4$.  The direct product of the Fibonacci substitution of Example \ref{nonconst1D} with itself is shown in Figure \ref{DP1}.   
\begin{figure}[ht]
\includegraphics[width = 4.25in]{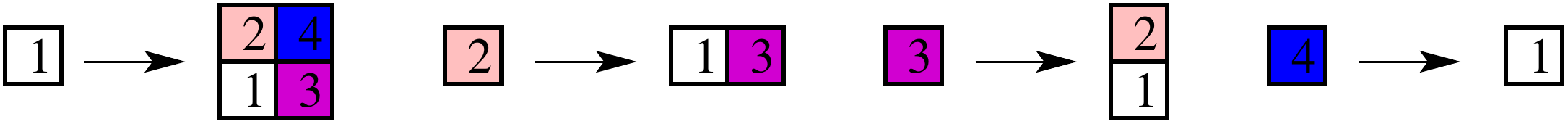}
\caption{The Fibonacci direct product substitution.}
\label{DP1} 
\end{figure}
Using only colors without the numbers we show the level-0 through level-4 blocks of type $1$ in Figure \ref{DP2}. 
The characteristic ``plaid" appearance of the direct product is evident.

\begin{figure}[ht]
\includegraphics[width = 4.75in]{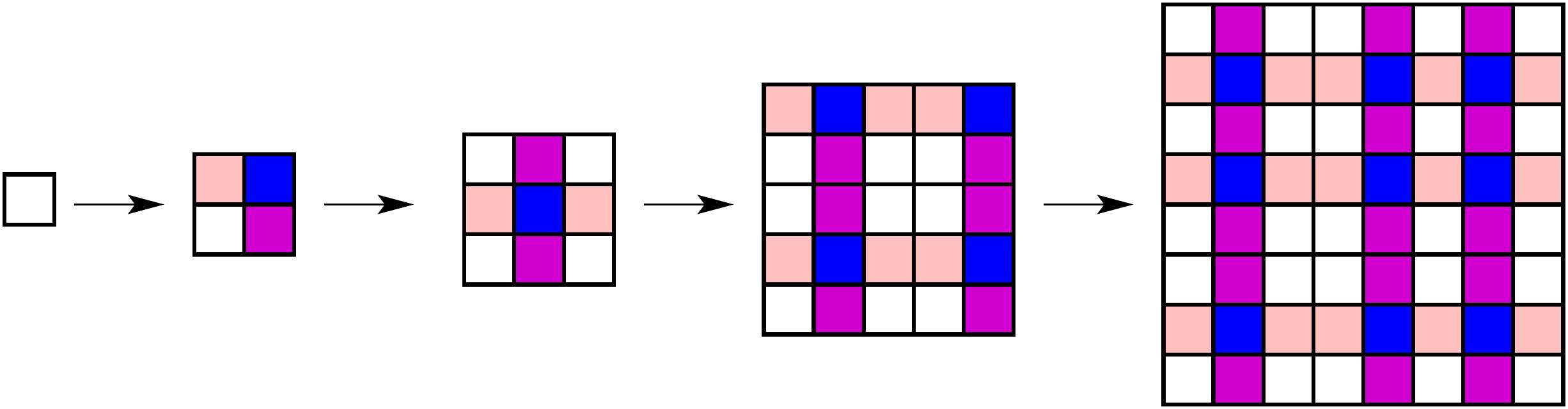}
\caption{A few iterations of the Fibonacci direct product substitution.}
\label{DP2} 
\end{figure}

\end{ex}
Some literature on $d$-dimensional symbolic substitutions exists.   In the non-constant length case, direct product substitutions, with a generalization allowing randomness in the choice of substitution from level to level, are studied in \cite{Mozes}.   An extension of this idea, allowing substitutions with restrictions forcing the substitutions to ``fit", are studied in \cite{Hansen}.   In the constant-length case, a partial survey and spectral analysis of this class from the dynamical systems viewpoint appears in \cite{Mytoppaper}.    For those wishing to experiment with various substitutions of both constant and non-constant length,  the author maintains a MATLAB freeware computer program that allows the user to generate these tilings of $\Z^2$ and manipulate them in several ways \cite{tiling.program}.

\subsection{Tilings of $\R^d$}  Let us introduce some terminology that will be useful throughout the paper.   A {\em tile} is a set $T \subset \R^d$ that is the closure of its interior.   We will always assume that tiles are bounded; in the literature it is frequently assumed that tiles are connected or even homeomorphic to topological balls.   In fact it is often required that the tiles be polygonal, but in substitution tiling theory tiles with fractal boundary occur naturally.   
When it is desirable to distinguish between congruent tiles they can be {\em labeled} (also called {\em marked} or {\em colored}).   Two tiles are considered {\em equivalent} if they differ by a rigid motion and carry the same label.   
A {\em prototile set} is a finite set $\ppp$ of inequivalent tiles.
Given a prototile set $\ppp$, a {\em tiling } of  $\R^d$ is a set $\T$ of tiles, each equivalent to a tile from $\ppp$,  such that
\begin{enumerate}
\item $\T$ covers $\R^d$:   $\R^d = \bigcup \{T : T \in \T\}$ , and
\item $\T$ packs $\R^d$: distinct tiles have non-intersecting interiors.  
\end{enumerate}

A {\em $\T$-patch} is a finite union of tiles  with nonintersecting interiors covering a connected set; two patches are equivalent if there is a rigid motion between them that matches up equivalent tiles.  
A tiling is said to be of {\em finite local complexity (FLC)} (also known as having a {\em finite number of local patterns}) if there are only finitely many two-tile $\T$-patches up to equivalence.   A tiling is called {\em repetitive} (also called {\em almost periodic} or the {\em local isomorphism property}) if for any $\T$-patch $P$ there is an $R>0$ such that in every ball of radius $R$ there is a patch equivalent to $P$.
In dynamical systems theory the most work has been done on repetitive tilings with finite local complexity.

\subsection{Infinite tilings from substitutions; tiling spaces and dynamical systems} Given a tiling substitution, it is possible to construct infinite tilings and tiling spaces from that substitution in a few different ways.   (This is also true for symbolic substitutions).   Our description will be necessarily imprecise as different substitutions can require  different definitions of some of the terms; we give the main ideas here and refer the reader to sources such as \cite{Radinmiles, Robinson.ams}, and \cite{Sol.self.similar} to get more details.  

One way to get an infinite tiling is to begin with some initial block or tile and substitute ad infinitum.   In many cases a limiting sequence or tiling $\T_0$ will exist.   Sometimes it will cover only a half-line, quarter-plane, or some other unbounded region of space, and sometimes it will cover the entire line or plane.
A less constructive method is to define a tiling $\T$ as {\em admitted} by the substitution if every finite configuration of tiles in $\T$ is equivalent to a configuration found inside a level-$n$ tile, for some $n$.

The {\em tiling space} associated to a substitution is the set of all tilings admitted by that substitution.   
Another way to obtain this space is to take the closure (in a suitable metric) of all rigid motions of a limiting tiling $\T_0$.
In either case,  a point in the tiling space $X$ is an infinite tiling, and any nontrivial rigid motion of that tiling is considered a different point in the tiling space.

\subsection{Outline of the paper}  Substitutions of constant length have a natural generalization to tilings in higher dimensions, which we introduce in Section \ref{geom.section}.     These generalizations, which include the well-studied self-similar tilings, rely upon the use of linear expansion maps and are therefore rigidly geometric.   We present examples in varying degrees of generality and include a selection of the major results in the field.

Extending substitutions of non-constant length to higher dimensions seems to be more difficult, and is the topic of Section \ref{comb.section}.   To even define what this class contains has been problematic and there is not yet a consensus on the subject.   For lack of existing terminology we have decided to call this type of substitution combinatorial as tiles are combined to create the substitutions without any geometric restriction save that they can be iterated without gaps or overlaps, and because in certain cases it is possible to define them in terms of their graph-theoretic structure.  

In many cases one can transform combinatorial tiling substitutions into geometric ones through a limit process.  In Section \ref{connections.section}, we will discuss how to do this and what the effects are to the extent that they are known.   We conclude the paper by discussing several of the different ways substitution tilings can be studied, and what sorts of questions are of interest.
\section{Geometric tiling substitutions}\label{geom.section}

Although the idea had been around for several years, self-similar tilings of the plane were given a formal definition and introduced to the wider public  by Thurston in a series of four AMS Colloquium lectures, with lecture notes appearing thereafter \cite{Thurston}.   
Throughout the literature one finds varying degrees of generality and some commonly used restrictions.   We make an effort to give precise definitions here, adding remarks which point out some of the differences in usage and in terminology.

\subsection{Self-similar tilings:  proper inflate-and-subdivide rules}
\label{sst.definitions}
For the moment we assume that the only rigid motions allowed for equivalence of tiles are translations;  this follows \cite{Thurston} and \cite{Sol.self.similar}.   We give the definitions as they appear in \cite{Sol.self.similar}, which includes that of \cite{Thurston} as a special case.

Let $\phi: \R^d \to \R^d$ be a linear transformation, diagonalizable over $\C$, that is expanding in the sense that all of its eigenvalues are greater than one in modulus.   A tiling $\T$ is called {\em $\phi$-subdividing} if
\begin{enumerate}
\item for each tile $T \in \T$, $\phi(T)$ is a union of $\T$-tiles, and
\item $T$ and $T'$ are equivalent tiles if and only if $\phi(T)$ and $\phi(T')$ form equivalent patches of tiles in $\T$.
\end{enumerate}

A tiling $\T$ will be called {\em self-affine with expansion map $\phi$} if it is $\phi$-subdividing, repetitive, and has finite local complexity.   If $\phi$ is a similarity the tiling will be called {\em self-similar}.   For self-similar tilings of $\R$ or $\R^2 \cong \C$ there is an {\em expansion constant} $\lambda$ for which $\phi(z) = \lambda z$.

The rule taking $T\in \T$ to the union of tiles in $\phi(T)$ is called an {\em inflate-and-subdivide rule} because it inflates using the expanding map $\phi$ and then decomposes the image into the union of tiles on the original scale.   If $\T$ is $\phi$-subdividing, then it will be invariant under this rule, therefore we show the inflate-and-subdivide rule rather than the tiling itself.   The rule given in Figure \ref{square1}   is an inflate-and-subdivide rule with 
$\phi(z) = 3z$.   However, the rule given in Figure \ref{DP1} is not an inflate-and-subdivide rule.

\begin{ex}
\label{chair}
The ``L-triomino" or ``chair" substitution uses four prototiles, each being an L formed by three unit squares.  We have chosen to color the prototiles since they are inequivalent up to translation.  The expansion map is $\phi(z) = 2z$ and in Figure \ref{Chair1} we show the substitution of the four prototiles.\begin{figure}[ht]
\includegraphics[width=3.5in]{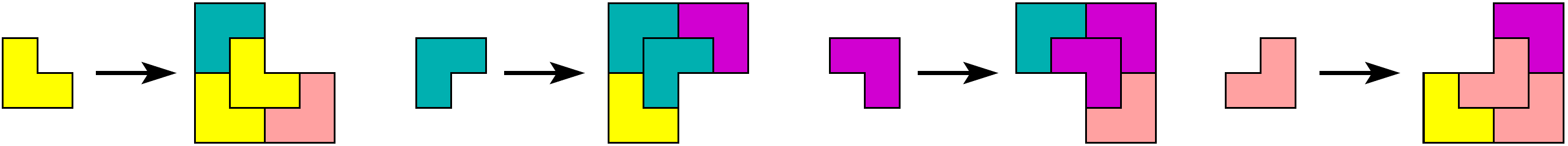}
\caption{The ``chair" or ``L-triomino" substitution.}
\label{Chair1} 
\end{figure}

This geometric substitution can be iterated simply by repeated application of $\phi$ followed by the appropriate subdivision.   Parallel to the symbolic case, we call a tile that has been inflated and subdivided $n$ times a {\em level-$n$ tile}.     In Figure \ref{chair2} we show level-$n$ tiles for $n = 2, 3,$ and $4$.

\begin{figure}[ht]
\includegraphics[width=5in]{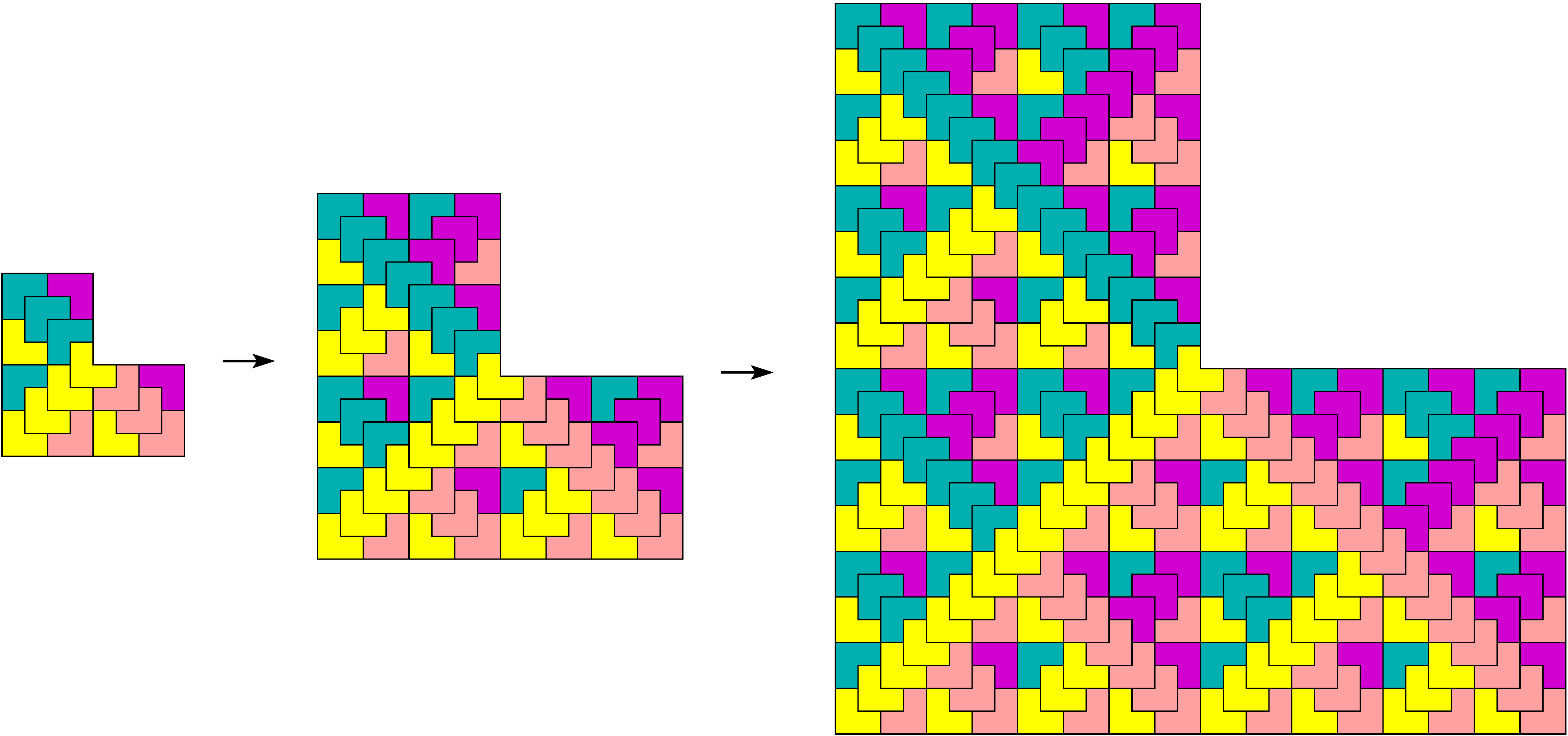}
\caption{ Level-2, level-3, and level-4 tiles.}
\label{chair2}
\end{figure}
\end{ex}

\subsection{A few important results} \label{sstresults}
One of the earliest results was a characterization of the expansion constant $\lambda \in \C$ of a self-similar tiling of $\C$.   
\begin{thm} (Thurston \cite{Thurston}, Kenyon \cite{Kenyon}) A complex number $\lambda$ is the expansion constant for some self-similar tiling if and only if $\lambda$ is an algebraic integer which is strictly larger than all its Galois conjugates other than its complex conjugate.
\end{thm}
The forward direction was proved by Thurston and the reverse direction by Kenyon.  In \cite{Kenyon2}, Kenyon extends the result to self-affine tilings of $\R^d$ in terms of eigenvalues of the expansion map.
%

In the study of substitutions, from one-dimensional symbolic substitutions to very general tiling substitutions, the {\em substitution matrix} is an indispensable tool.    (This matrix has also been called the ``transition", ``composition", ``subdivision", or even ``abelianization" matrix).  Suppose that the prototile set (or alphabet) has $m$ elements labeled by $\{1, 2, ..., m\}$
.   The {\em substitution matrix} $M$ is the $m \times m$ matrix with entries given by 
\begin{equation}
M_{ij} = \text{ the number of tiles of type } i \text{ in the substitution of the tile of type } j.
\end{equation}

For example, the substitution in Example \ref{const2D} has substitution matrix $M = \left(\begin{array}{cc} 9&1 \\ 0&8\end{array}\right)$ when we label $a = 1$ and $b = 2$.  If an initial configuration of tiles has $n$ white tiles and $m$ blue tiles, then $M [n \,\,m]^T$ is the number of white and blue tiles after one application of the substitution.

Since the substitution matrix is always an integer matrix with nonnegative entries, Perron-Frobenius theory is relevant (see for example \cite{Kitchens, Sol.self.similar}).  The results we need require $M$ to be {\em irreducible}:  for every $i, j \in \{1, 2, ..., m\}$ there exists an $n$ such that $(M^n)_{ij}>0$.   Among other things, the Perron-Frobenius theorem states that if $M$ is irreducible, then the largest eigenvalue will be a positive real number that is larger in modulus than any of the other eigenvalues of the matrix.   This eigenvalue is unique, has multiplicity one, and is called the {\em Perron eigenvalue} of the matrix.

Primitivity, a special case of irreducibility, is particularly important.  A matrix $M$ is {\em primitive} if there is an $n>0$
such that $M^n$ has strictly positive entries.  Primitivity of $M$ means if one substitutes any tile (or letter) a fixed number of times, one will see all of the other tiles (or letters).   This is a relatively strong property, and one that is almost always assumed in the literature.   All of the substitutions in the paper are primitive except Examples \ref{const2D} and \ref{Chacon1}.  

The following theorem, mentioned in \cite{Kenyon2}, is proved in \cite{Praggastis}.
\begin{thm}(Praggastis \cite{Praggastis})
An FLC $\phi$-subdividing tiling is repetitive if and only if its substitution matrix is primitive.
\end{thm}

Solomyak's papers \cite{Sol.self.similar} and \cite{Sol.u.comp} give several other key results for the dynamical systems of self-similar or self-affine tilings.  
The following is stated as a corollary to the Perron-Frobenius theorem.
\begin{thm} (Solomyak \cite{Sol.self.similar}) If $M$ is primitive, the Perron eigenvalue of $M$ is the volume expansion $|\text{det}\,\,\phi|$. The Perron left eigenvector gives the relative volumes of the prototiles.
\end{thm}
 
An algebraic integer is a {\em Pisot number} if all of its algebraic conjugates (except its complex conjugate) are smaller than one in modulus.   Whether or not the expansion constant $\lambda$ is a Pisot number is especially important from a dynamical point of view.   In dimensions one and two, Solomyak \cite{Sol.self.similar} has shown that a self-similar tiling dynamical system is not weakly mixing if and only if its expansion constant is a complex Pisot number.  We will see this number-theoretic property having other effects in Sections \ref{comb.section} and \ref{connections.section}.

So far, the results in this section have depended only on the substitution matrix and expansion constant, and not the geometry of the substitution.  The final theorem in this section uses the notion of ``matching rules", which are fundamentally geometric.   Roughly speaking, a set of matching rules determine which patches are allowed in a tiling.  A simple yet classic example is the Penrose tiling with marked rhombs, which we will encounter in Example \ref{Penrose.rhombs}.    The markings give matching rules that ``enforce" the substitution in the sense that any tiling of the plane constructed following the matching rules must be admitted by the Penrose substitution.    Goodman-Strauss was able to generalize this result to most geometric tiling substitutions:
\begin{thm}(Goodman-Strauss \cite{Goodman-Strauss}) Every (geometric) substitution tiling of $\R^d$, $d > 1$ can be enforced with finite matching rules, subject to a mild condition:   We require that tiles admit a set of ``hereditary edges" such that the substitution tiling is ``sibling edge-to-edge".
\end{thm}
We leave a discussion of the particulars to \cite{Goodman-Strauss} and note only that the ``sibling edge-to-edge" condition is mild enough to encompass most of the known examples.
\subsection{Geometric generalization:  infinite rotations and sizes}
In the previous sections the only rigid motions allowed for tile equivalence were translations.    However,
there are natural tiling substitutions that require relaxing this to allow rotations.   

\begin{ex}
\label{pinwheel}
The ``pinwheel" substitution rule acts on right triangles of
side lengths 1, 2, and $\sqrt{5}$, inflating them by a factor of $\sqrt{5}$ and subdividing into
five triangles as shown in Figure \ref{pinwheel1}.

\begin{figure}[ht]
\includegraphics[width=3.5in]{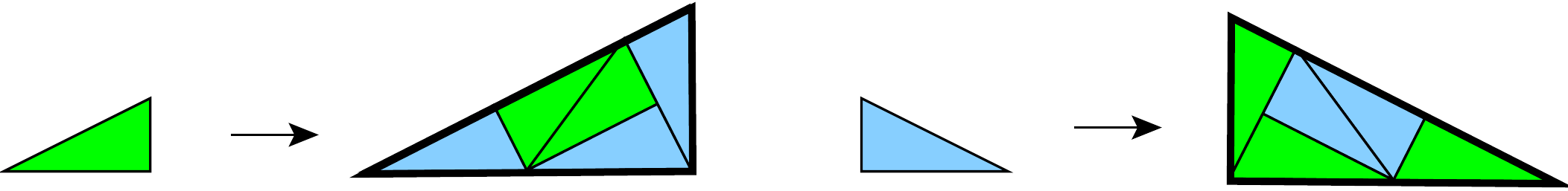}
\caption{The pinwheel substitution}
\label{pinwheel1} 
\end{figure}

Radin introduces pinwheel tilings in \cite{Radinpinwheel}, attributing them to unpublished work of John H. Conway.  Radin proves that the tiles appear in an infinite number of orientations
that are uniformly distributed mod $2\pi$, and calls the tiling space ``statistically round" because it is invariant under rotations of infinite order.
He also establishes matching rules for the pinwheel substitution, giving us
the first amazing example where local matching rules produce infinite rotations!

\begin{figure}[ht]
\includegraphics[width=3.25in]{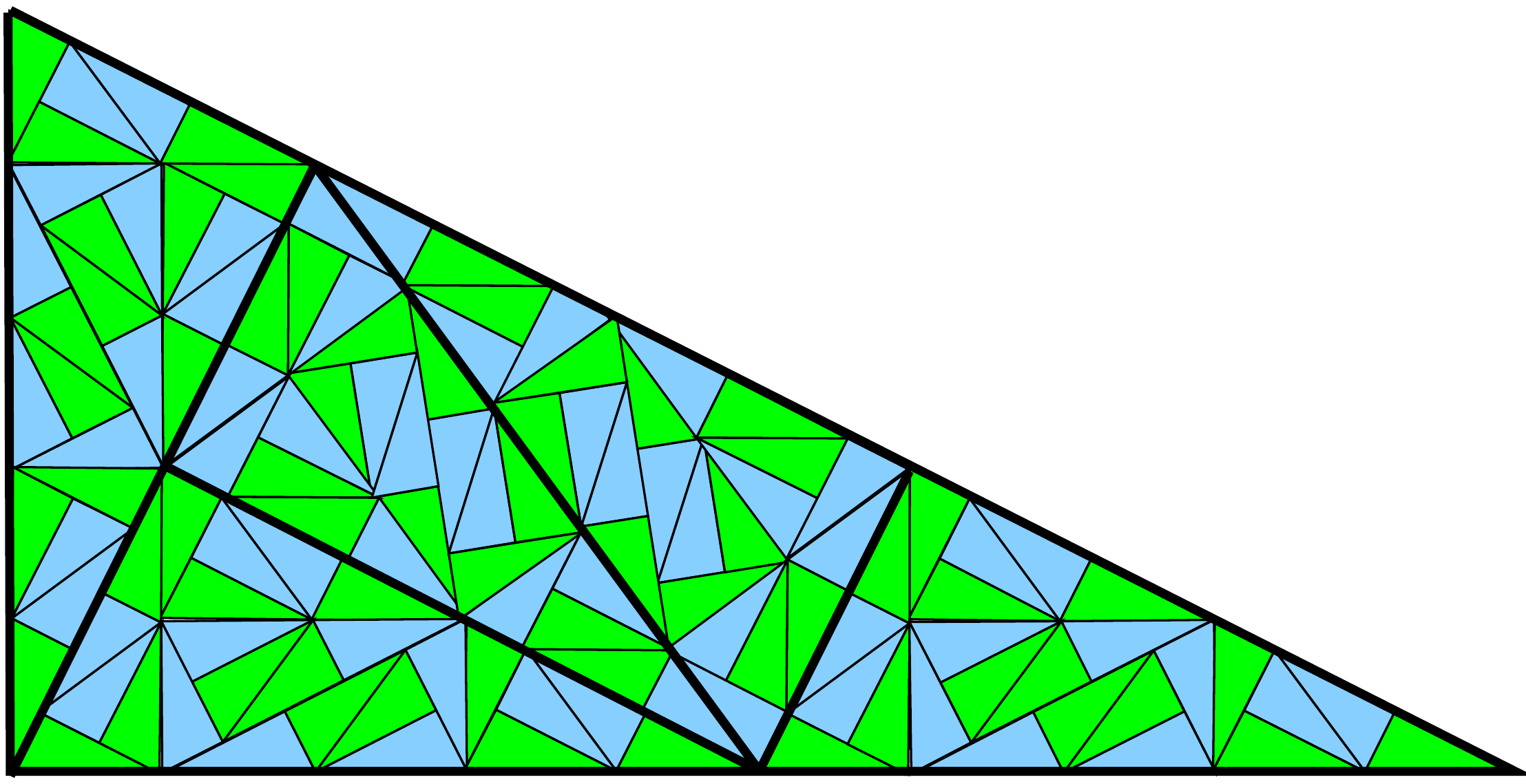}
\caption{A level-3 tile for the pinwheel inflation.}
\label{pinwheel2} 
\end{figure}

\end{ex}

\begin{ex}\label{pinwheelgenex}
Sadun \cite{Sadun.pin} comes up with an interesting twist on the previous example: the ``generalized pinwheel" tilings.  Instead
of requiring that the tiles be isometric to members of some finite prototile set, he requires only that
they be conformally equivalent.   The subdivision rules are quite straightforward but some choices arise for the inflation portion and we do not attempt to explain those here.  At the first level, one takes a right triangle with side lengths $a, b$ and $c$ and subdivides it into 5 similar triangles as in Figure \ref{pinwheelgen}. 
\begin{figure}[ht]
\includegraphics[width=2in]{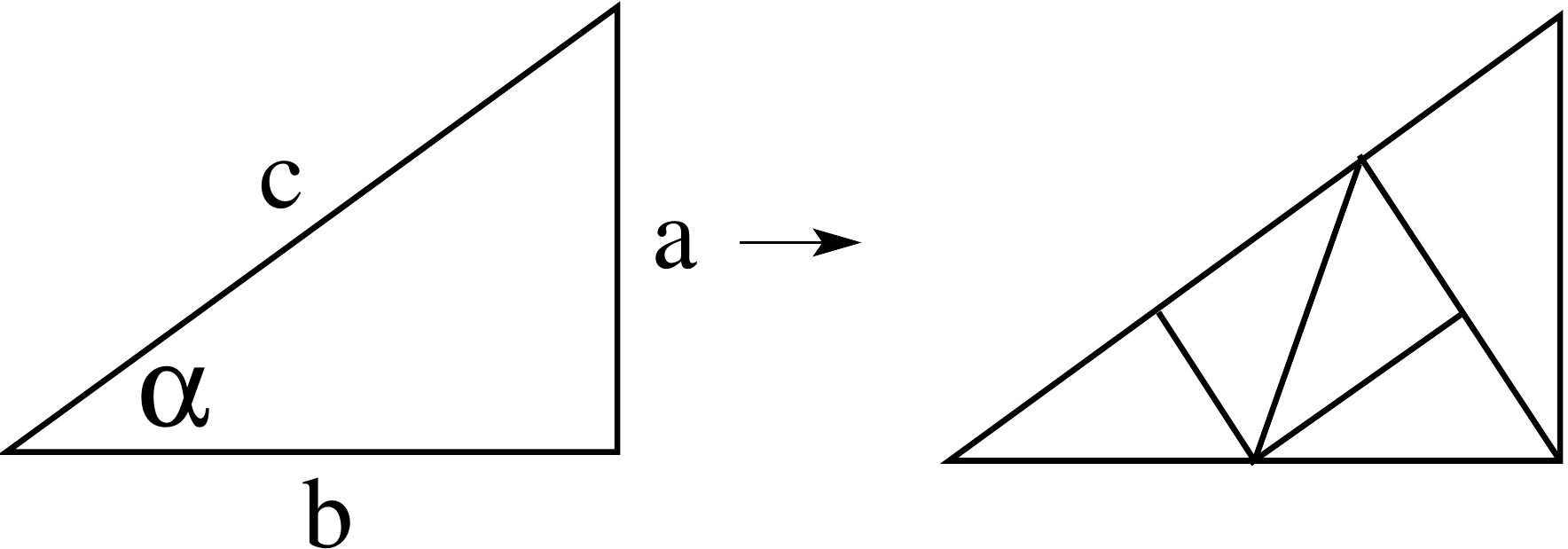}
\caption{The first decomposition in Sadun's pinwheel generalization.}
\label{pinwheelgen} 
\end{figure}
The subdivision at the next stage takes place only on the largest of the triangles.  One can continue subdividing indefinitely; for the appropriate inflation at any stage and for precision of the results mentioned below we refer the reader to \cite{Sadun.pin}.

All salient properties of a tiling admitted by the substitution depend on the angle $\alpha$.
For the original pinwheel tiling $\alpha = \sin^{-1}(1/\sqrt{5})$, and that is one of the angles for which the
tiling has finite local complexity, but the tiles appear in an infinite number of orientations.  For other values
of $\alpha$ the tiles will only appear in a finite number of sizes.  There is only one value, $\alpha = \pi/4$, for which the tiles appear in both a finite number of sizes and a finite number of orientations.   Figure \ref{pinwheel4} shows a few subdivisions of this case.

\begin{figure}[ht]
\includegraphics[height=.7in]{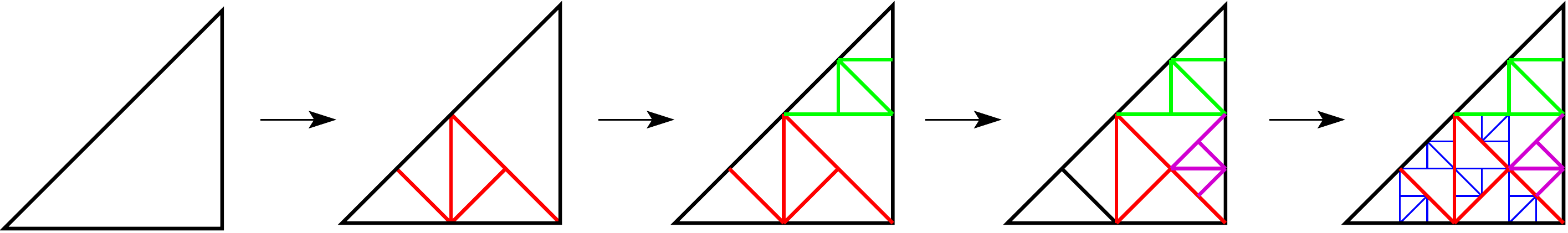}
\caption{A special case of the generalized pinwheel tilings: $\alpha = \pi/4$.}
\label{pinwheel4}
\end{figure}
\end{ex}

\subsection{Geometric generalization: pseudo-self-similar tilings}

A close cousin of the self-similar tiling is the pseudo-self-similar tiling, which is generated by a variant of the inflate-and-subdivide rule.  Tiles are still inflated and then replaced by tiles from the original scale, but these may stick out of or not completely cover the inflated tile, so the substitution rule is ``imperfect".   We will show two well-known substitutions in this section,  both of which can be converted into proper inflate-and-subdivide rules in different ways.   

\begin{ex} \label{Penrose.rhombs}The Penrose inflation using marked rhombs is shown in Figure \ref{Penroseinfl}.   When the reader attempts to inflate and subdivide a second time, she will notice that the subdivisions of adjacent tiles overlap.   This is not a contradiction, however, because the overlapping tiles are equal and will therefore be considered the same tile.

\begin{figure}[ht]
\includegraphics[height=1in]{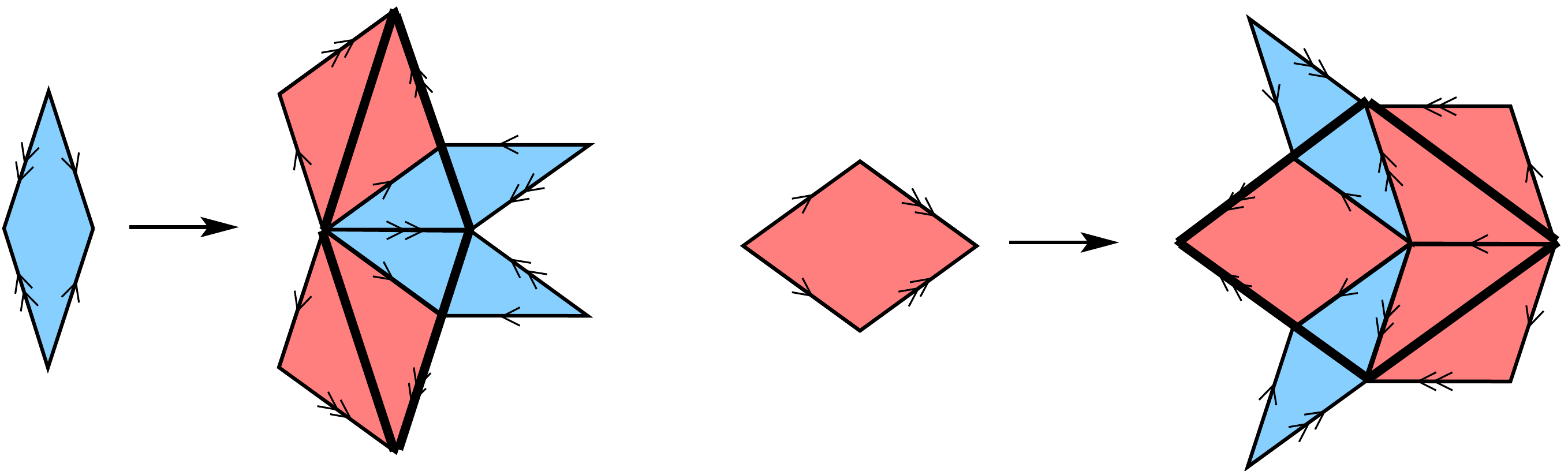}
\caption{The ``Penrose inflation": an imperfect substitution rule.}
\label{Penroseinfl}
\end{figure}

\begin{figure}[ht]
\includegraphics[width=6in]{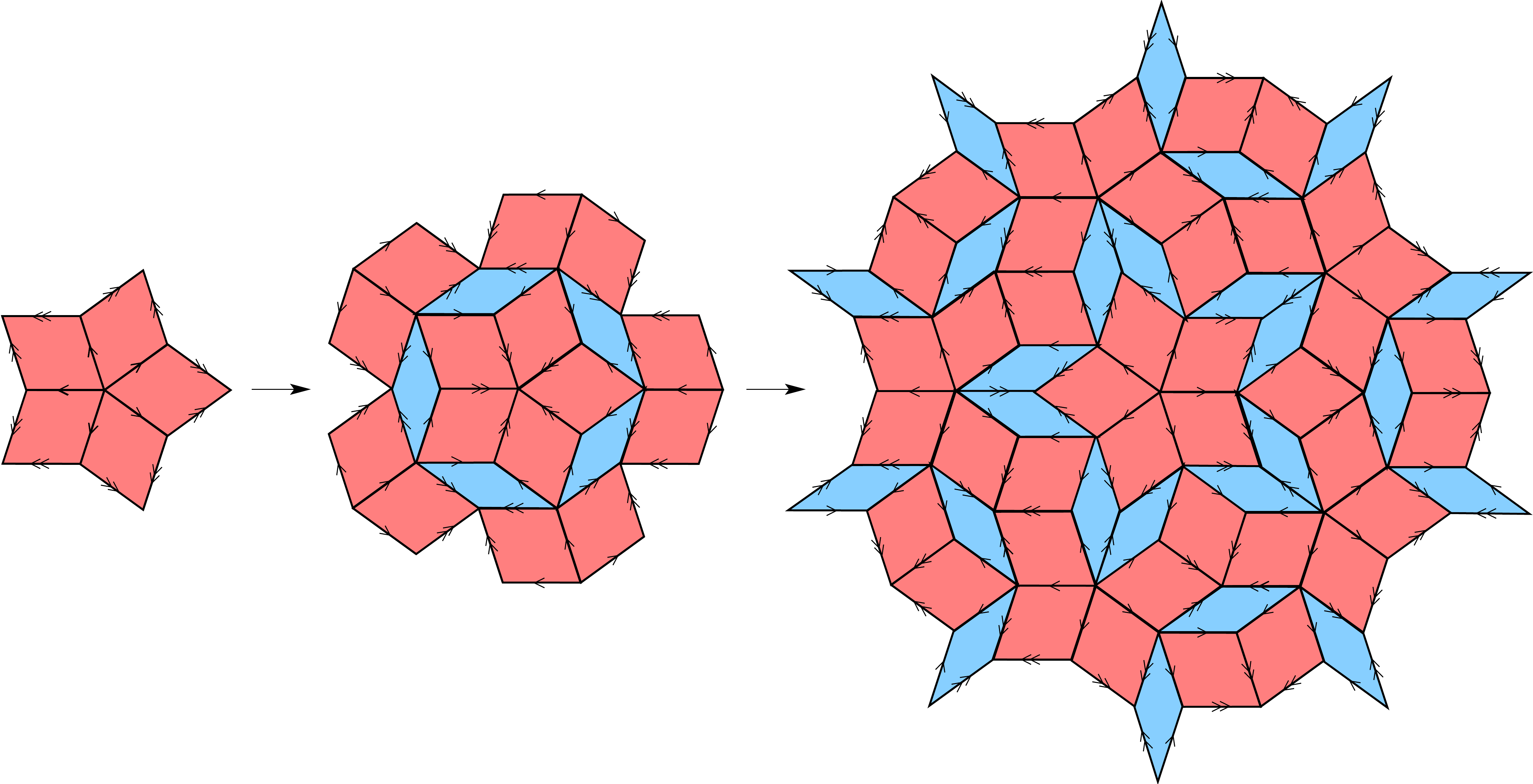}
\caption{A few iterations of the Penrose inflation.}
\label{Penroseits}
\end{figure}

Penrose tilings appear in many equivalent forms, with alternative tile shapes such as triangles or the famous ``kites and darts".    The Scientific American article by Gardner \cite{Gardner}
introduced the Penrose tilings and many of their interesting properties to the general public.
 In Chapter 10  of \cite{G-S}, the Penrose tiles are studied as kites and darts and in other forms. 
 It is possible to slice the Penrose rhombs in half, creating triangles (known as Robinson triangles) on which a proper inflate-and-subdivide rule can be defined; we show this in Figure \ref{Penrose3}.  Chapter 6 of \cite{Sen} gives a detailed analysis of Penrose tilings and includes the Robinson triangles.

 \begin{figure}[ht]
\includegraphics[height=1.25in]{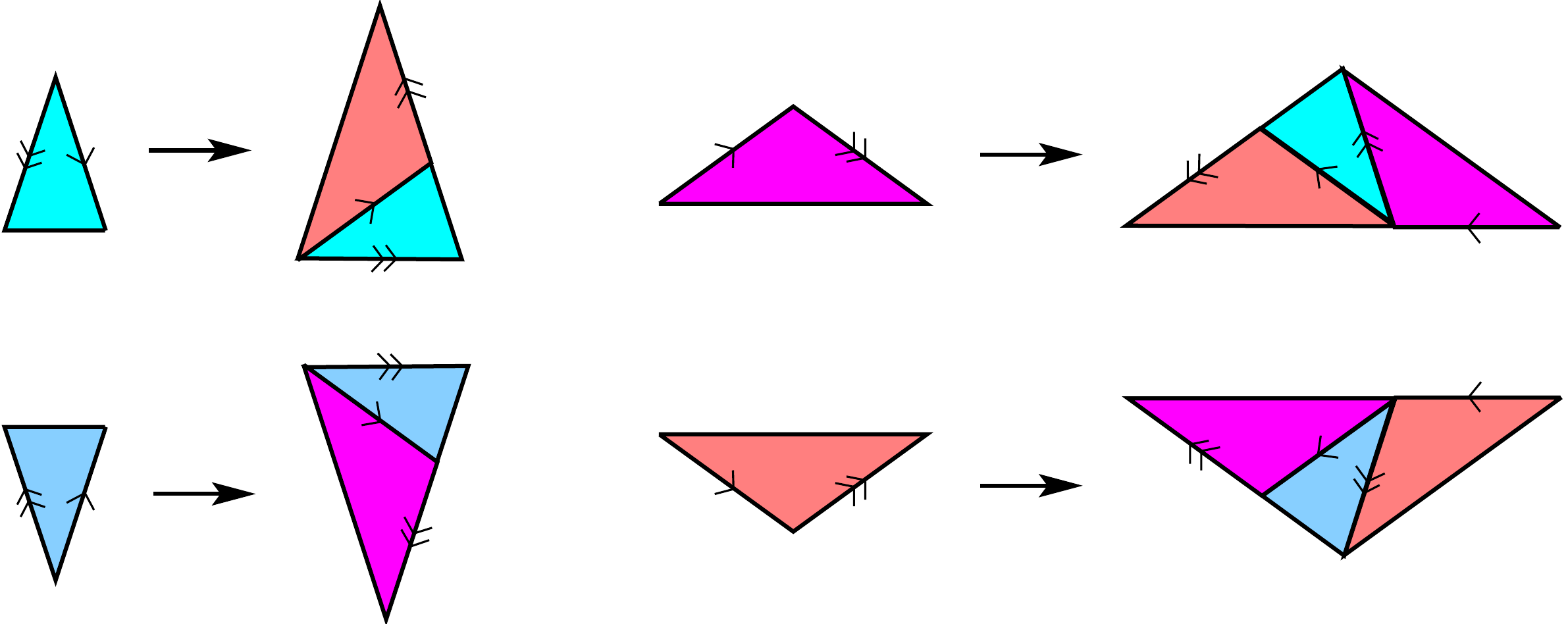}
\caption{The triangles version of the inflate-and-subdivide rule.}
\label{Penrose3}
\end{figure}

\end{ex}

Penrose tilings have a number of interesting properties, most of which can be found in other tilings but were first observed in the Penrose tilings.   We mention a few of the highlights and leave the details to the references.  
There are (at least!) three ways, other than substitution, to generate Penrose tilings.   One, which we have already discussed, is via matching rules: if an infinite tiling is constructed from Penrose tiles with the requirement that adjacent tiles have matching arrows (both in number and direction), this tiling will be a Penrose tiling.    Amazingly, this local activity of matching the arrows forces the global property of being generated by a substitution!  

Two other methods for generating Penrose tilings are the multigrid method and the projection method.   Since both methods rely on lattices, they can be used to prove that the nonperiodicity of the Penrose tilings is a tightly controlled form of disorder.
The {\em multigrid method} was discovered by DeBruijn \cite{DeBruijn}.  
 In this method one superimposes five grids of lines to create a ``pentagrid"; 
every pentagrid is dual to a Penrose tiling.   This method is used by E. A. Robinson, Jr. \cite{Rob.penrose} to understand the Penrose tilings as a dynamical system.
A nice description of the {\em projection method} appears in \cite{Sen}, p. 195-196.   In this method a copy of $\R^d$ is embedded in $\R^{d+n}$,  and some lattice $\lll \subset \R^{d+n}$ is chosen.   Points from $\lll$ are projected onto the copy of $\R^d$ to form a tiling of $\R^d$.    Projection tilings have been
studied extensively and are only sometimes obtainable by methods of substitution.   A characterization of the intersection is given in a certain case in \cite{Harriss-Lamb}.  Point sets generated by generalized projection methods are called ``model sets" and are of great interest in mathematics and physics, and their spectrum is the subject of intense study (see e.g. \cite{Lee.Moody.Solomyak,Moody.survey}).

If the reader were to experiment with a large set of Penrose rhombs, he would quickly discover that it is difficult to tile a large region without disobeying the matching rules.  Indeed, almost all finite configurations of Penrose rhombs, no matter how large, cannot be extended to infinite tilings \cite{Penrose.incompletable}!  Thus the fact that an infinite Penrose tiling exists at all is a major result, and it can be proved using the presence of the inflation rule.

\begin{ex}\label{Binaryex} The ``binary tilings" (see \cite{Sen}) are generated by the substitution rule shown in Figure \ref{Binary1} using the unmarked Penrose tiles.  The volumes expand linearly by a factor of $(5+\sqrt{5})/2$. This substitution is interesting from the dynamical viewpoint as it produces a weakly mixing tiling dynamical system.  The fact that it is weakly mixing means that it has a level of disorder not present in the Penrose tilings, despite the fact that the relative frequencies of thick to thin rhombs is the same in both tilings!  Weak mixing is evident in the diffraction spectrum, which is pictured in \cite{Sen} and is analyzed in a few papers including
\cite{Godreche}.  \begin{figure}[ht]
\parbox{.55in}{
\includegraphics[width=.5in]{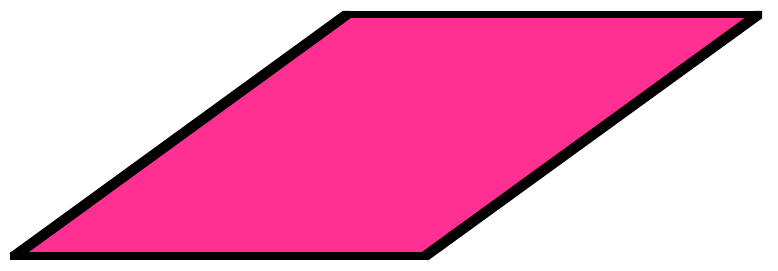}}
{\includegraphics[width=.25in]{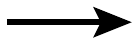}}
\parbox{1.5in}{\includegraphics[width=1.15in]{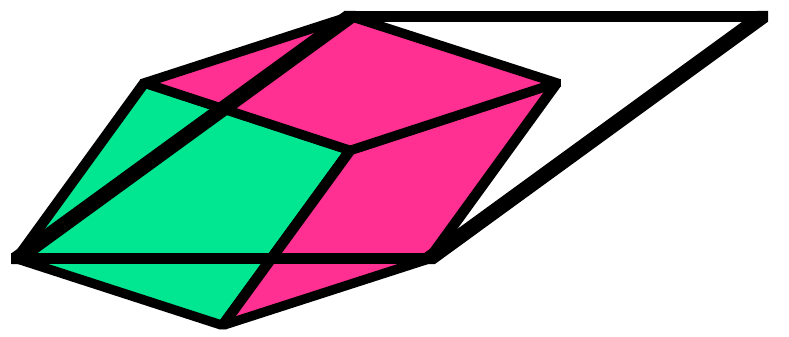} } \hspace{1cm} 
\parbox{.45in}{\includegraphics[width=.4in]{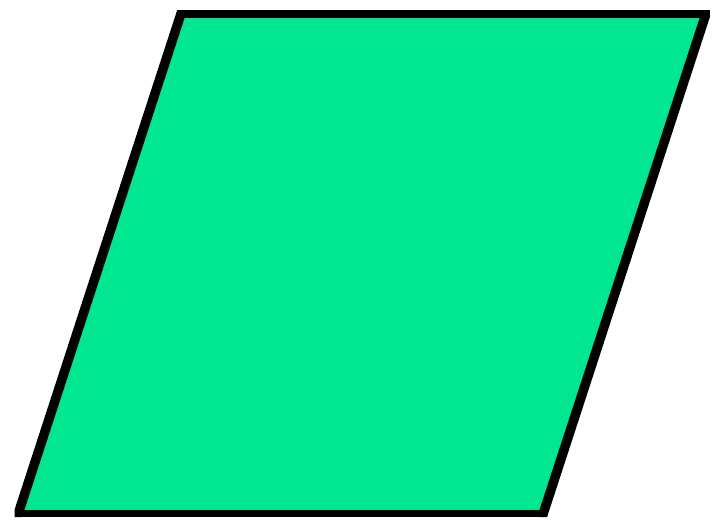}}
\includegraphics[width=.25in]{arrow}
\parbox{1.15in}{\includegraphics[width=.85in]{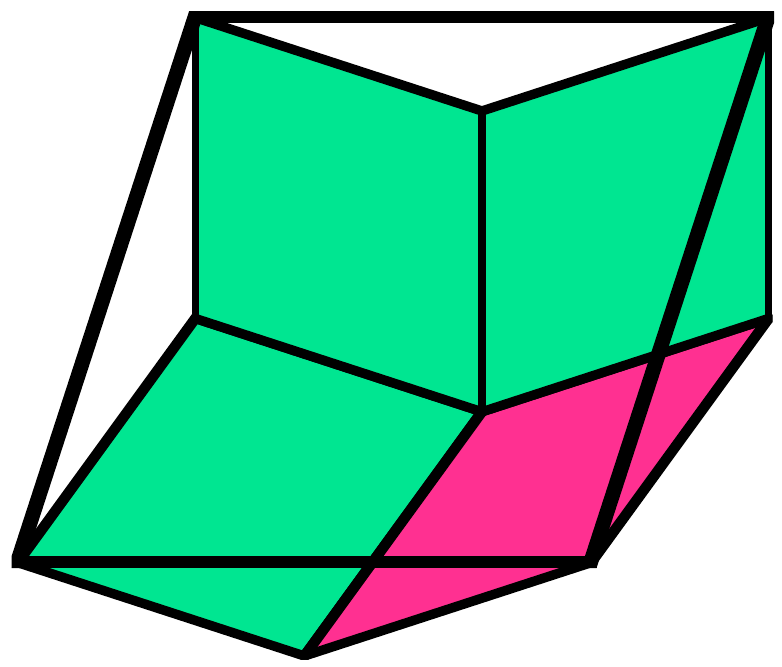}}
\caption{Another imperfect substitution uses the unmarked Penrose rhombs.}
\label{Binary1} 
\end{figure}

\begin{figure}[ht]
\parbox{.55in}{\includegraphics[width=.55in]{binarysub1}}
\includegraphics[width=.25in]{arrow}
\parbox{1.25in}{\includegraphics[width=1.25in]{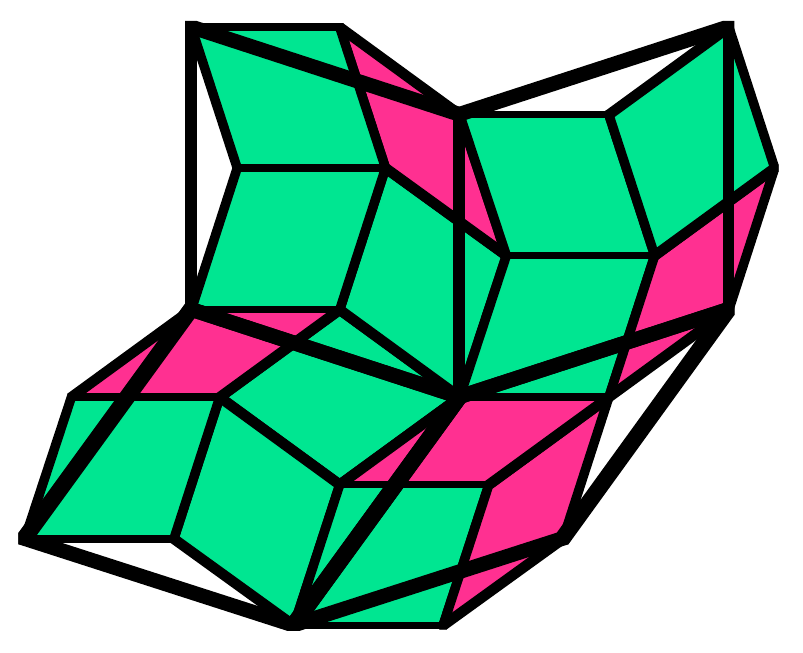}}
\includegraphics[width=.25in]{arrow}
\parbox{2.5in}{\includegraphics[width=2.5in]{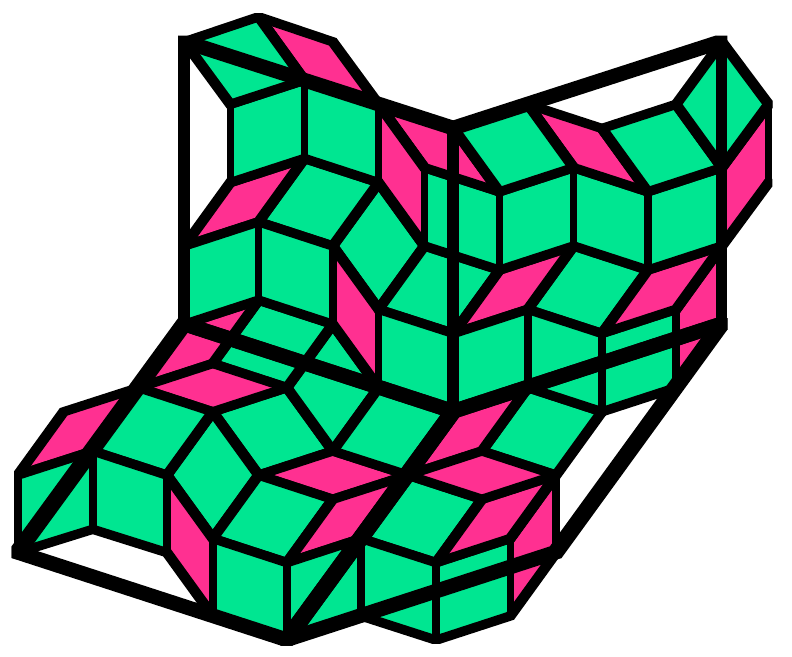}}
\caption{How to iterate the binary substitution.}
\label{Binary2} 
\end{figure}
\end{ex}

\begin{figure}[ht]
\includegraphics[width=5in]{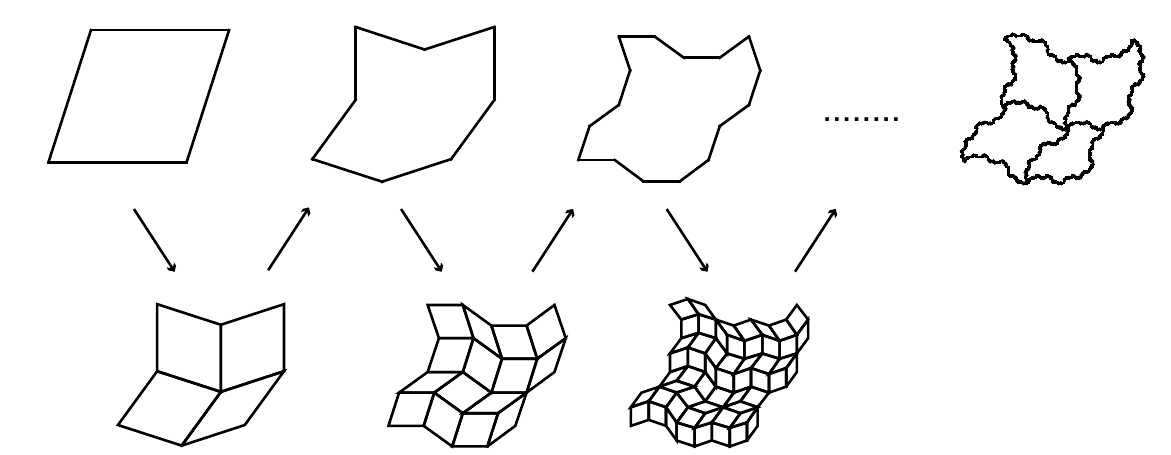}
\caption{Obtaining fractal binary tiles.}
\label{Binaryfractal}
\end{figure}

One can ``redraw" the tiles of the binary tiling in such a way that the substitution rule becomes
a perfect inflate-and-subdivide rule that gives rise to a self-similar tiling of the plane.    The method
is to take the limit  as $n$ goes to infinity of the support of the level-$n$ tiles divided by the $n$th power of the expansion constant (see Figure \ref{Binaryfractal}, courtesy of E. Arthur Robinson, Jr.).
A general version of this process has been shown to work for all pseudo-self-similar tilings
of the plane \cite{Meboris}.   This result is extended to tilings of $\R^d$ in \cite{Sol.psst}, provided the tiles appear in a finite number of orientations.  In Section \ref{connections.section} we will use a similar process to transform combinatorial substitutions into geometric ones, with some unexpected results.

 \section{Combinatorial tiling substitutions}
 \label{comb.section}
 
The self-similar tilings and their close relatives in the previous section come from substitution rules that have one thing in common:  a single similarity (or expanding linear map, in the self-affine case) governs the inflation of all of the tiles.  Now we consider substitutions for which a tile and its replacement may be geometrically unrelated, or for which there are several linear maps governing the tile inflations.   There is no unified definition for this class of substitutions.  Attempts to define a tiling substitution based on the dual graph of the tiling have been made \cite{Peyriere,Mycstpaper, ABS}, and we call this method ``constructive".  However, there are perfectly reasonable tiling substitutions for which combinatorial information is insufficient to define the substitution rules. Substitutions of this type we call ``non-constructive".   Leaving a formal presentation of the definitions to the references, we simply present ideas and some examples.

\subsection{Constructive combinatorial tiling substitutions}
In this type of substitution, level-$n$ tiles can be constructed using only the information about adjacencies between tiles, making it possible to iterate the substitution using only local information.

\begin{figure}[ht]
\includegraphics[height=2cm]{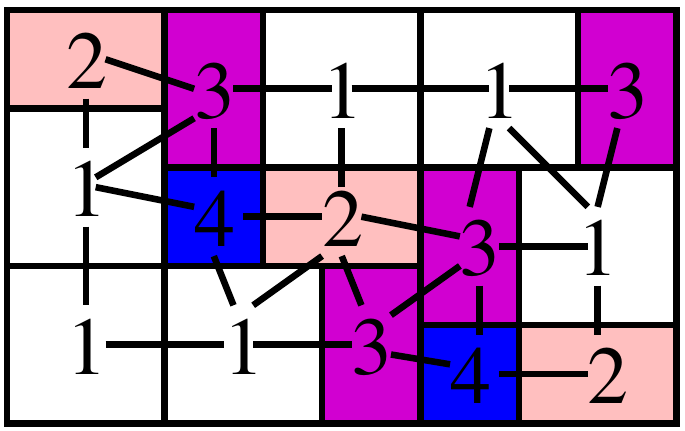}
\caption{Dually situated tiling and graph.}
\label{dual.graph}
\end{figure}

A tiling of $\R^2$ makes a drawing of a planar graph:  there is a vertex wherever three or more tiles meet, an edge wherever exactly two tiles meet, and the tiles themselves are the facets.  The elements of this graph can be labeled according to the tile and adjacency types they represent; we may choose a labeling scheme that provides as much or as little information about the surrounding tiles as we wish.   Any planar tiling has a {\em dual graph}: vertices of the dual graph are tiles of the tiling, there is an edge between two vertices if the corresponding tiles are adjacent, and facets correspong to vertices in the tiling.   The elements of a dual graph inherit the labels of their dual counterparts.  Part of a tiling and its dual graph are shown in Figure \ref{dual.graph}, using the conventions that numbers represent both the label and the vertex of the graph, and that the edge and facet labels are suppressed.

The dual graph is a natural object on which to define a substitution, by analogy with one-dimensional symbolic substitutions.  The labeled vertices are the letters of an alphabet that corresponds to our prototile set. Any rule that replaces a vertex with a finite graph might correspond to a substitution rule for a prototile, provided the geometry of the tiles allows it.  A difficulty is that two vertices (tiles) may be adjacent in many different ways, so to make the dual graph ``see" this, we keep track of the adjacency types in our edge and facet labels.  Constructive combinatorial substitutions specify exactly how the substituted graphs of two adjacent elements should be attached.  The next example, taken from \cite{Mycstpaper}, will provide some intuition.

\begin{ex}\label{fib.dpv} We obtain a {\em direct product variation} (DPV) substitution by rearranging some of the tiles in the Fibonacci substitution of Figure \ref{DP1} to break up the direct product structure.  Here we have carefully rearranged the tiles in the substitution of the type 1 tile so that the substitution can be iterated without inconsistency.   We show the result in Figure \ref{graph.sub}, along with the induced graph substitution. 

\begin{figure}[ht]
\includegraphics[width=4in]{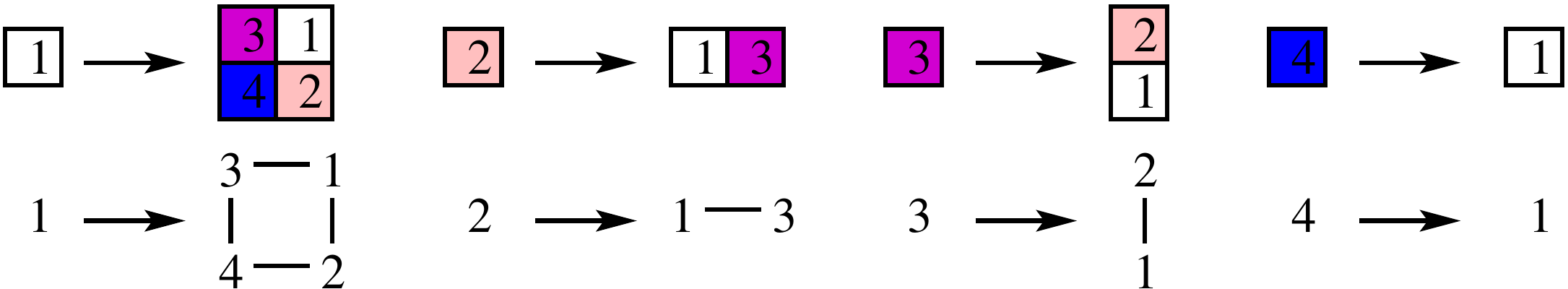}
\caption{The vertices of the dual graph inherit the substitution.}
\label{graph.sub}
\end{figure}

To obtain the level-$n$ block, one simply concatenates the level-$(n-1)$ blocks in the ``obvious" manner as shown in Figure \ref{DPV2}, matching sides that have the same length in the order prescribed by Figure \ref{graph.sub}.  Since the side lengths of the level-$(n-1)$ blocks are Fibonacci numbers, the fit is guaranteed at each stage.   Note that the ``plaid" appearance of the direct product (Figure \ref{DP2}) has disappeared.
\begin{figure}[ht]
\includegraphics[width = 4.5in]{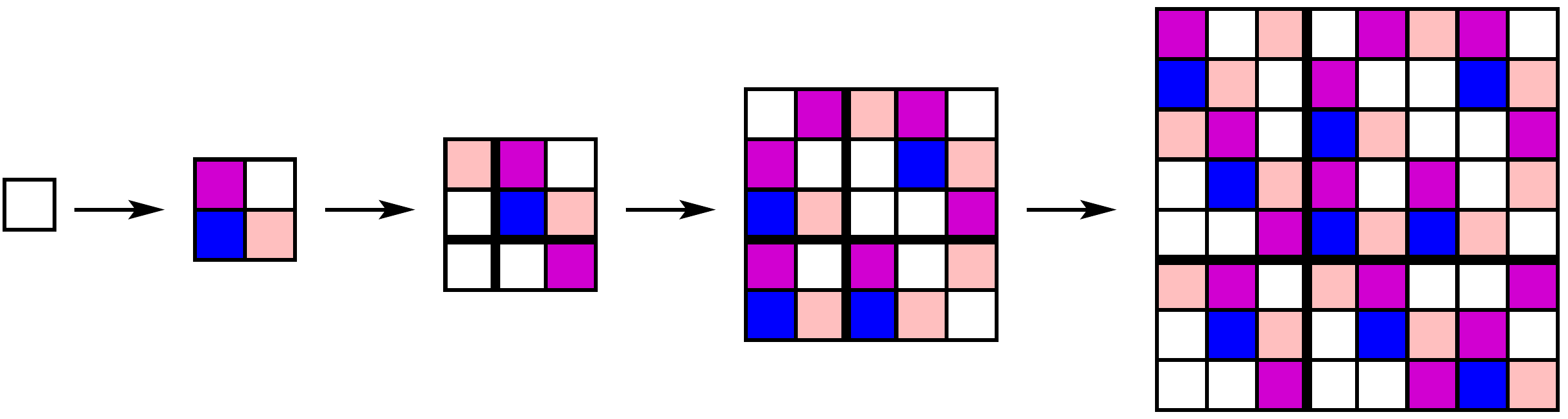}
\caption{Iterating of the Fibonacci DPV on the tile of type 1.}
\label{DPV2} 
\end{figure}

But what if you want to substitute a pair of adjacent tiles within a level-$n$ block?   You cannot do it consistently without knowing the larger context of the adjacency, that is, the tiles that surround them in the tiling.   For example, consider two horizontally adjacent tiles of type 1.   That adjacency appears twice in the level-3 tile of type 1 and we have circled them on the left side of Figure \ref{adj.subs}; what happens under substitution is shown on the right of the same figure.  For this graph substitution, the problem can be handled by relabeling the edges and facets of the graph in terms of the immediate configuration of tiles the edge or facet is contained in.

\begin{figure}[ht]
\includegraphics[width = 2.75in]{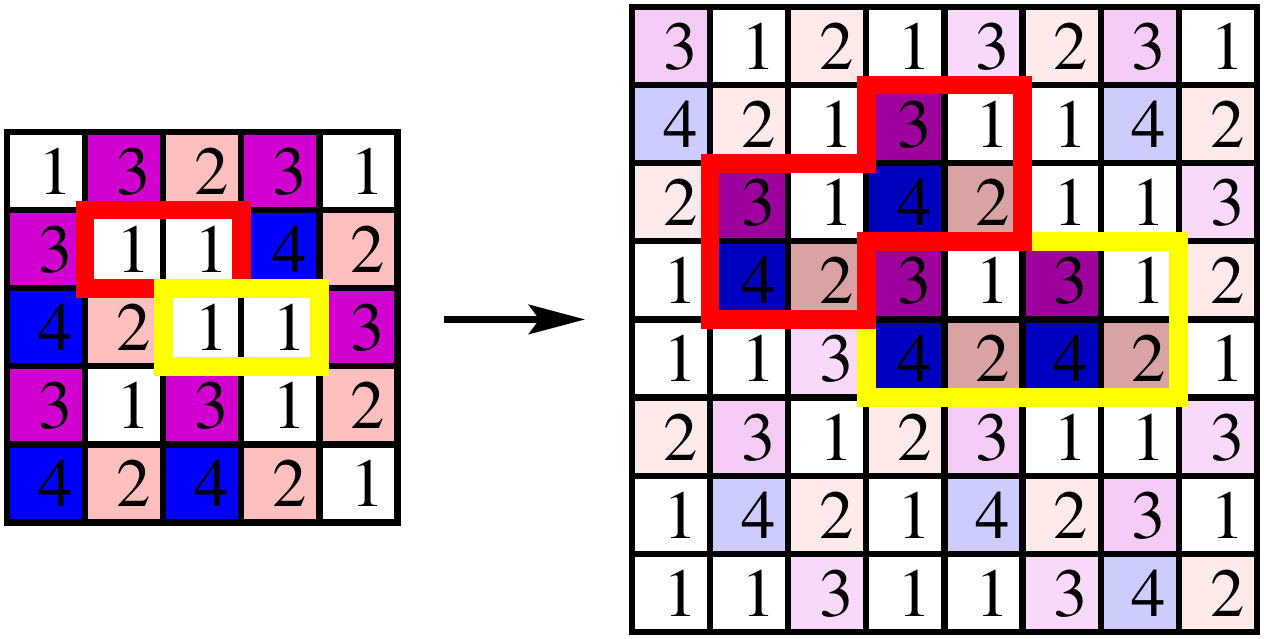}
\caption{The substitution of an adjacent pair of tiles depends on its context.}
\label{adj.subs}
\end{figure}

\end{ex}

The basic idea of a constructive combinatorial substitution for tilings as it appears in \cite{Mycstpaper} and in a similar form in \cite{ABS} is this.   Given a labeled vertex set $V$ representing the prototile types,  a map from $V$ to the set of nonempty labeled graphs on $V$ is the basis for the substitution rule.   The edges and facets of these graphs are labeled to give information about the types of adjacencies they represent in a tiling.   The substitution rule also specifies how to substitute the labeled edges and facets so that we know how to connect the vertex substitutions contained in certain labeled graphs.   (The need to specify graph substitutions on facets and not just edges is illustrated in an example in \cite{Mycstpaper}).  Defining a tiling substitution rule this way is quite tricky since most labeled graphs do not represent the dual graph of a tiling.    This interplay between combinatorics and geometry is where the technicalities come in to the formal definitions in the literature.

\begin{ex}\label{Rauzy.cst}
The tiling substitution of Figure \ref{Rauzy1}, introduced in \cite{ABS}, is based on a variation of the one-dimensional ``Rauzy substitution" $\sigma(1) = 1\,2, \,\,\sigma (2) = 3$, $\sigma(3) = 1$.  \begin{figure}[ht]
\includegraphics[width=2.5in]{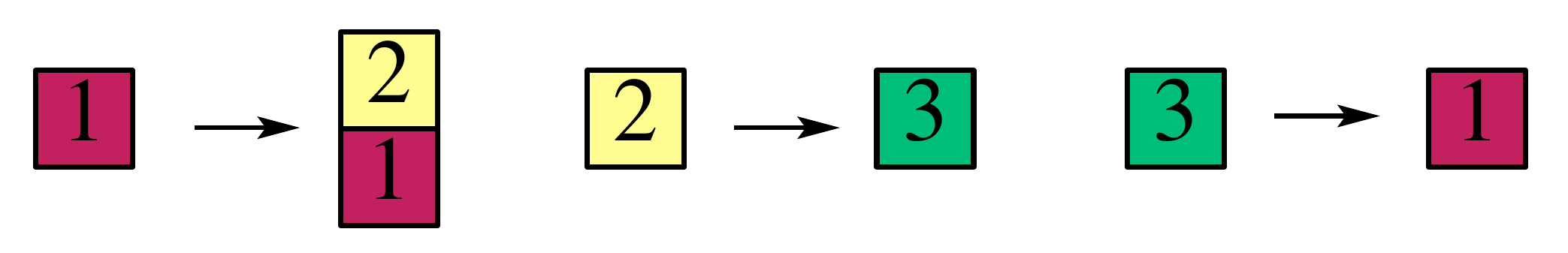}
\caption{A two-dimensional substitution based on the Rauzy one-dimensional substitution.}
\label{Rauzy1} 
\end{figure}
Figure \ref{Rauzy1} is obviously not enough information to iterate the substitution, so we specify how to substitute the ``important adjacencies" in Figure \ref{Rauzy2}.   This is enough \cite{ABS}: there are no ambiguities when substituting other adjacencies, and facet substitutions do not include any new information.  We show a few iterates of the tile of type 1 in Figure \ref{Rauzy3}, starting with the level-2 tile of type 1.
\begin{figure}[ht]
\includegraphics[width=3.75in]{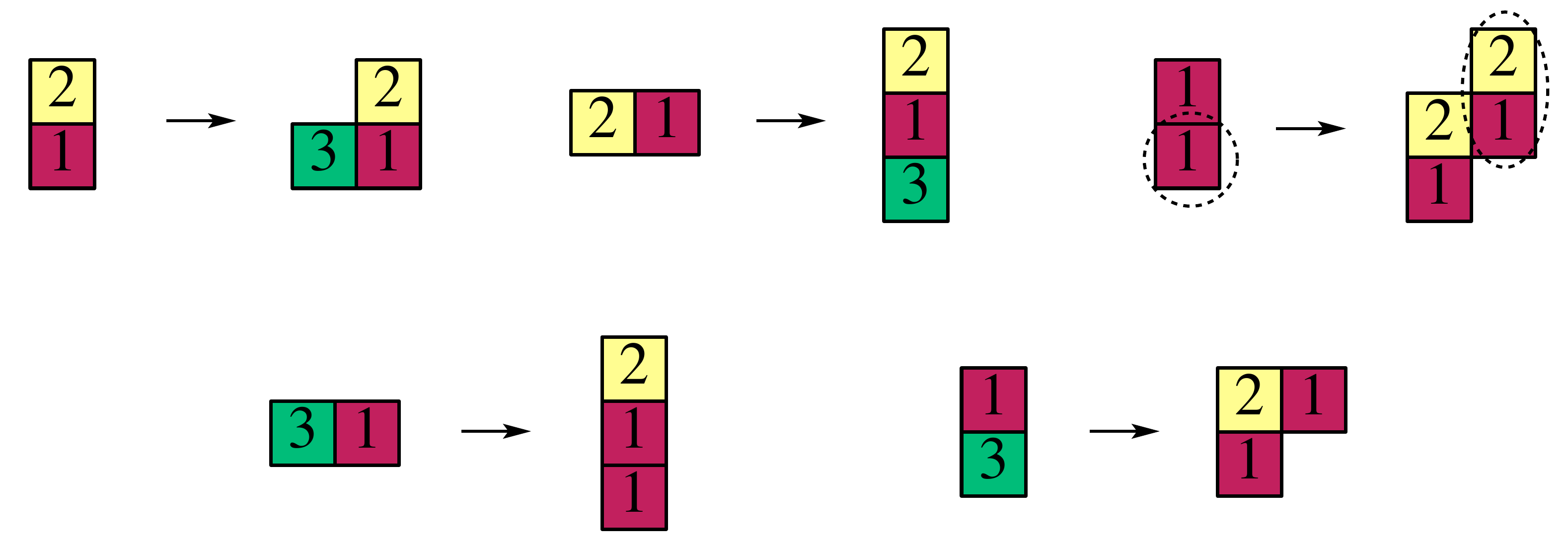}
\caption{How to substitute important adjacencies for the Rauzy substitution.}
\label{Rauzy2} 
\end{figure}
The fact that this substitution rule can be extended to an infinite tiling of the plane is proved
using noncombinatorial methods in \cite{ABS}; a combinatorial proof of existence would be welcomed.
\begin{figure}[ht]
\includegraphics[width=6in]{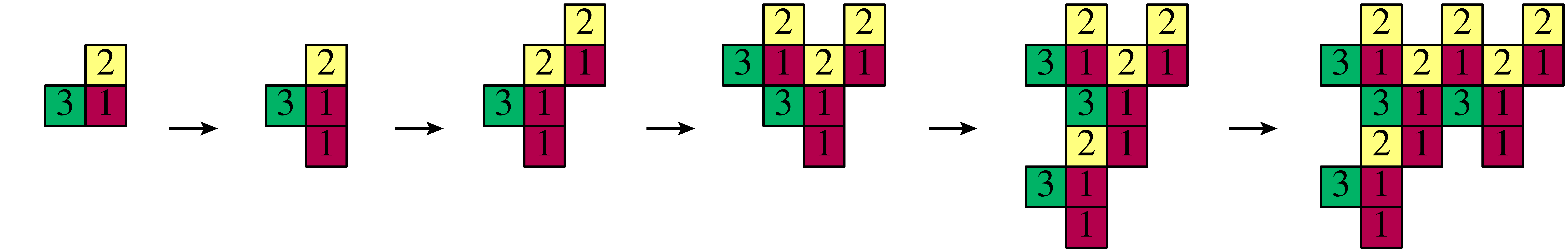}
\caption{A few iterates of the Rauzy two-dimensional substitution.}
\label{Rauzy3} 
\end{figure}
\end{ex}

\subsection{Non-constructive tiling substitutions}
When trying to make up new examples of combinatorial tiling substitutions it is easy to create examples that fail to be constructive.  The problem arises in the substitution of adjacencies: it may happen that no finite label set can be chosen to describe all adjacencies  sufficiently to know how to substitute them.   There is evidence to suggest that this sort of example can arise when the constant which best approximates the linear growth of blocks is not a Pisot number.  
The author is not aware of any formal definition containing this group and so proposes the following definition, which works directly with the tiling and does not involve dual graphs.
\begin{dfn}   \label{non.constructive.cst} A {\em (non-constructive) tiling substitution} on a finite prototile set $\ppp$ is a set of nonempty, connected patches $\sss=\{S^n(p) : p \in \ppp \text{ and } n \in 1,2,...\}$ satisfying the following:
\begin{enumerate}
\item  For each prototile $p \in \ppp$ and tile $t \in S^1(p)$, and for each integer $n \in 2, 3, ...$, there are rigid motions $g(p,n,t):\R^d \to \R^d$ such that
$ S^n(p)=\bigcup_{t \in S^1(p)} g(p,n,t)\left(S^{n-1}(t)\right)$, where
\item for any $t \neq t' $ in $S^1(p)$, the patches $g(p,n,t)\left(S^{n-1}(t)\right)$ and $g(p,n,t')\left(S^{n-1}(t')\right)$ intersect at most along their boundaries.
\end{enumerate}
\end{dfn}
We say a tiling $\T$ is {\em admitted by the substitution $\sss$} if every patch in $\T$ appears as a subpatch of some element of $\sss$.
This very general definition is satisfied by every substitution appearing in this paper except the Penrose substitution on rhombi and those generalized pinwheel tilings that do not allow a finite number of tile sizes.

The Rauzy substitution of Example \ref{Rauzy.cst} has a particularly efficient representation by this definition.  The patches $S^1(p)$ are given to the right of the arrows in Figure \ref{Rauzy1}, with all lower right corners at the origin.  Now, $S^n(2)=S^{n-1}(3)$ and $S^n(3) = S^{n-1}(1)$, so $g(2,n,3)$ and $g(3,n,1)$ are the identity map.   We find $S^n(1) = S^{n-1}(1) \cup g(1,n,2)\left(S^{n-1}(2)\right)$, so $g(1,n,1)$ is the identity map and all we have left to figure out is the formula for $g(1,n,2)$.   It turns out that $g(1,n,2)$ is  translation by a vector $\vecv_n$ that can be computed recursively.  Let $\vecv_0=(0,0), \vecv_1=(0,1)$ and $\vecv_2 = (-1,0)$; for $n \ge 3$ we have that $\vecv_n = \vecv_{n-3} - \vecv_{n-2}$.

The Fibonacci DPV of Example \ref{fib.dpv} also has a relatively simple formuation in terms of Definition \ref{non.constructive.cst}.   The side lengths of level-$n$ tiles are given recursively, and the placement of the level-$(n-1)$ tiles to create level-$n$ tiles depends only on these side lengths.   Thus the translations $g(p,n,t)$ are computable recursively as well.
The next example is also a DPV, but it cannot be defined in terms of dual graphs and is non-constructive.   We encourage the reader to think about how to write up Definition \ref{non.constructive.cst} in this case.    

\begin{ex}\label{np.dpv}
Consider a DPV arising from a 
one-dimensional substitution $a \to a b b b, b \to a$.  
From the direct product of this substitution with itself, we choose only to rearrange the substitution of the type-1 tile as in Figure \ref{Non-Pisot1}. 
\begin{figure}[ht]
\includegraphics[height=.65in]{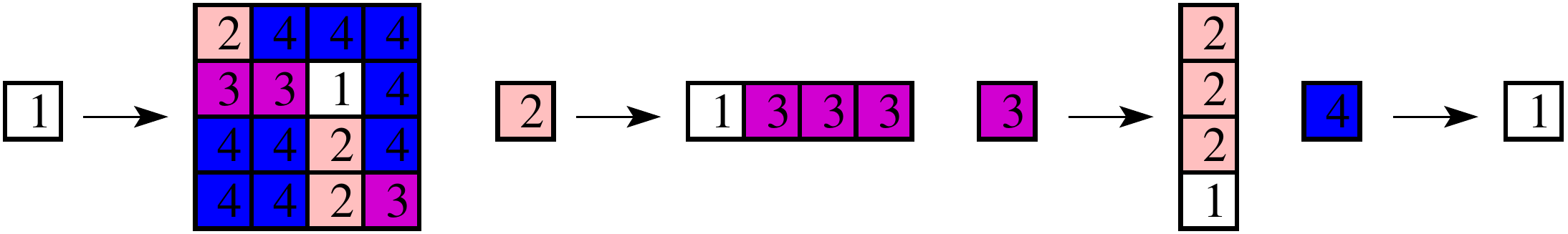}
\caption{A strongly non-Pisot direct product variation substitution.}
\label{Non-Pisot1} 
\end{figure}

The substitution matrix of this one-dimensional substitution has Perron eigenvalue $\theta=(1 + \sqrt{13}/2)$, which is not a Pisot number: its algebraic conjugate $\theta_2=(1 - \sqrt{13}/2)$ is larger than one in modulus.  Using analysis similar to what we will see in Section \ref{connections1D}, we find that constants times $\theta^n$ are the best approximation to the lengths of level-$n$ blocks.   The reader can jump ahead to the left side of Figure \ref{Non-Pisot3} to see the effect of $\theta$ being non-Pisot on the substitution.   There we find an adjacent pair of tiles (circled in red) which, under substitution, go to disjoint level-1 tiles (circled in green).  
It is possible to show that for any large $R>0$ there is a level-$n$ tile containing adjacent tiles $t_1$ and $t_2$ such that the distance between the substitutions of $t_1$ and $t_2$ in the level-$(n+1)$ tile is greater than $R$.  The method of proof relies on the fact that $\theta$ is non-Pisot and the result in \cite{Merobbie}.

\end{ex}

To illustrate that non-constructive tiling substitutions can have non-square tiles, 
to show an interesting connection to ``fault lines" (see Section \ref{curiousconnections}), and for its entertainment value, we include an unpublished substitution discovered by the author in 2002.    

\begin{ex}\label{nptriangles}
This substitution uses eight prototiles, shown to the left of the arrows  in Figure \ref{Non-pisottriangles1}.
\begin{figure}[ht]
\includegraphics[width=5in]{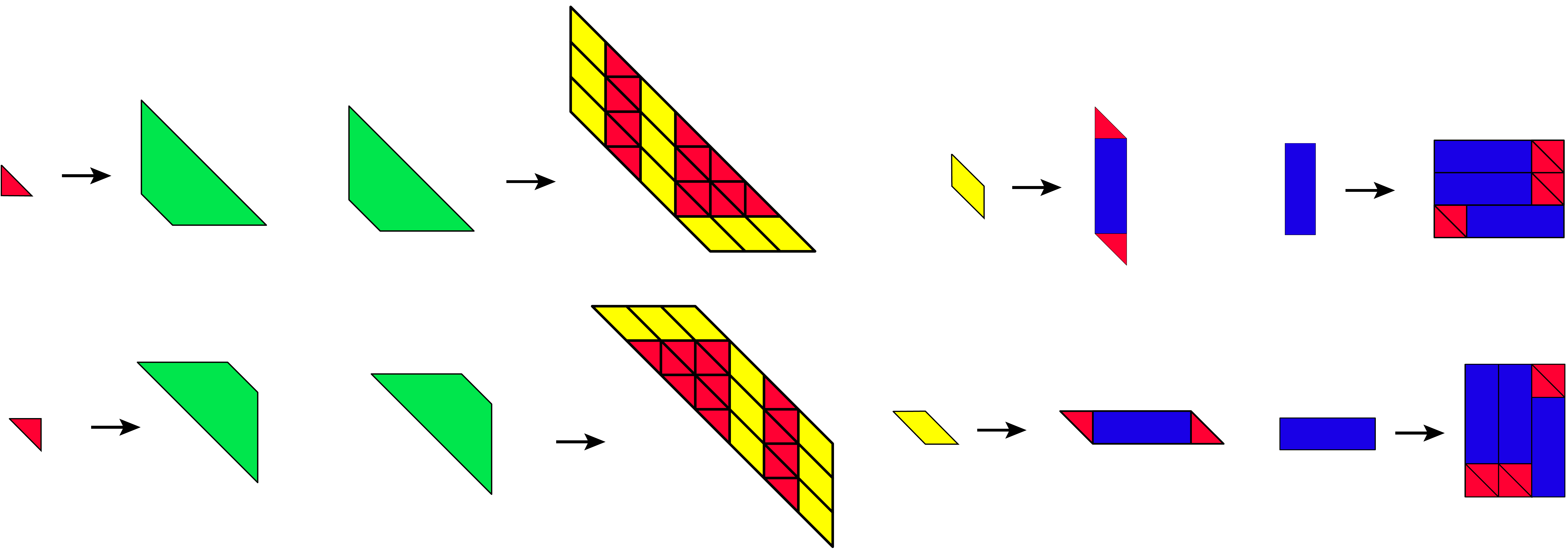}
\caption{A substitution similar to one in \cite{Mycstpaper}.}
\label{Non-pisottriangles1} 
\end{figure}
The Perron eigenvalue of the substitution matrix is the same non-Pisot number as in the previous example.  
Again an argument can be made to show that there are always adjacencies
that pull apart under the substitution, no matter how well we try to label them.   The left side of Figure \ref{Non-pisottriangles2} shows four iterations of the square formed by two red triangles.

 \end{ex}
 
 \section{Connections between geometric and combinatorial substitutions}
 \label{connections.section}
In Example \ref{Binaryex}, we saw an improper substitution that could be ``fixed" by redrawing the tile boundaries with the replace-and-rescale method.   In this method, the new tile boundary is the limit of the boundaries of the level-$n$ tiles, rescaled by the $n$th power of the inflation factor.   In Example \ref{Binaryex}, the result was a set of fractal prototiles satisfying a proper inflate-and-subdivide rule.   It turns out that this process can be used to create prototile sets for self-similar tilings in the case of some combinatorial substitutions as well.   We begin with the one-dimensional case, which is well-understood.
 
 \subsection{One-dimensional case}\label{connections1D}
 Given a symbolic substitution of non-constant length, 
it is easy to create an inflate-and-subdivide rule on labeled intervals (tiles) that has the same combinatorics as the symbolic substitution. The expansion constant is the Perron eigenvalue of the substitution matrix and the tile lengths are given by the Perron eigenvector.   We illustrate with the Fibonacci substitution in a way that looks ahead to the higher-dimensional case.
 
 \begin{ex}\label{nonconst1Dsst}
We introduce some geometry by thinking of two unit length prototiles represented by different colors in Figure \ref{Fib.subs.unit}.   \begin{figure}[ht]
 \includegraphics[width=3in]{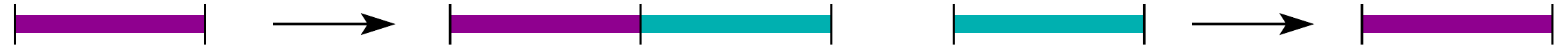}
 \caption{Fibonacci substitution with unit length tiles}
 \label{Fib.subs.unit}
 \end{figure}
The reader can check that the level-$n$ tiles will have lengths given by the entries of $(1,1)M^n$, where $M$ is the substitution matrix  $\left[\begin{array}{cc}1 & 1\\1&0\end{array}\right]$.   This matrix has eigenvalues given by $\frac{1 \pm \sqrt{5}}{2}$: the golden mean 
$\gamma$ and its conjugate $1-\gamma = -1/\gamma$.   Let $\vecv_1=(\gamma, 1)$ and $\vecv_2=(1, -\gamma)$ denote the associated left eigenvectors.   There are constants $k_1$ and $k_2$ so that $(1,1) = k_1\vecv_1 + k_2\vecv_2$, which gives us the vector of level-$n$ tile lengths:
\begin{equation} \label{Fib.length.eqn}
(1,1)M^n= \gamma^n k_1 \vecv_1 + (-1/\gamma)^n k_2 \vecv_2
\end{equation}
The lengths of the intervals for our self-similar tiling are the entries of $\lim_{n \to \infty} \gamma^{-n} (1,1)M^n= k_1 \vecv_1$.   The length of the type-$a$ tile is $k_1 \gamma$ and the length of the type-$b$ tile is $k_1$.   These lengths form an eigenvector for $M$, so
there exists an inflate-and-subdivide rule, which we have shown in Figure \ref{Fib.subs}.
 \begin{figure}[ht]
 \includegraphics[width=3in]{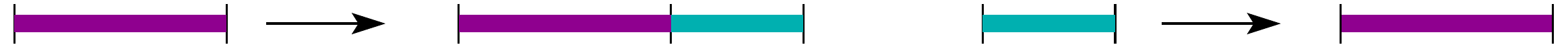}
 \caption{Fibonacci inflate-and-subdivide rule}
 \label{Fib.subs}
 \end{figure}
 \end{ex}
 
\noindent{\bf Notes:} (1) This process works on any substitution on $m$ letters provided that the vector $(1,1,...,1)$ lies in the span of the left eigenvectors of the substitution matrix of the substitution.  It works trivially on constant length substitutions since  $(1,1, ..., 1)$, the vector of unit tile lengths, already forms a Perron eigenvector for the substitution matrix.

(2) In the Fibonacci example, since $\gamma$ is a Pisot number (its conjugate $-1/\gamma$ is smaller than one in modulus), the higher the inflation the less important the second term of Equation 
 \ref{Fib.length.eqn} becomes.   Thus the lengths of the level-$n$ tiles of the inflate-and-subdivide rule are asymptotically close to the lengths of the non-constant length substitution, and therefore approximately integers!  The situation is dramatically different, of course, if any of the secondary eigenvalues are strictly greater than one in modulus (the strongly non-Pisot case).
 
 \subsection{Two-dimensional case}\label{connections2D}
 The reader should not be too surprised to discover that this process will work for direct product substitutions, such as Example \ref{nonconst2D}, and their variations, such as Examples \ref{fib.dpv} and \ref{np.dpv}.   The level-$n$ blocks are rectangular and have side lengths given by the lengths of the one-dimensional substitution.   Rescaling by the expansion factor gives us rectangular tiles whose side lengths are determined by the Perron eigenvector  as before.    Thus, if the original one-dimensional substitutions have a proper inflate-and-subdivide rule,  so will any DPVs associated with them.
Still, it is instructive to consider the two-dimensional replace-and-rescale method as it applies in this simple case.
 \begin{ex}\label{fib.sst}
Consider the Fibonacci DPV substitution of Example \ref{fib.dpv}.   There are four tile types and the substitution matrix is $M = \left[\begin{array}{cccc}1&1&1&1\\1&0&1&0\\1&1&0&0\\1&0&0&0\end{array}\right].$  If there is a proper inflate-and-subdivide rule corresponding to our substitution, it must have the same substitution matrix.  This tells us that the volume expansion of the rule must be the Perron eigenvalue, which in this case is $\gamma^2$, the square of the golden mean.   (The other eigenvalues are 1, 1, and $1/\gamma^2$.)
The level-$n$ tiles of the DPV substitution are supported on rectangles with side lengths given by either the $n$th or the $(n-1)$st Fibonacci numbers.  We rescale the volumes by $1/\gamma^{2n}$ to obtain prototiles for our self-similar tiling.   

The ``right" way to see this process is to consider a linear map $\phi: \R^2 \to \R^2$ that expands with the Perron eigenvalue of the substitution matrix $M$.   (How to find this map in general is quite unclear.) In this example $\phi$ is given by the matrix $\left[\begin{array}{cc}\gamma&0\\0&\gamma \end{array}\right]$.   Denoting the support of the level-$n$ tile of type $t$ as $supp(S^n(t))$, we can find the support of the prototile $t'$ for the inflate-and-subdivide rule that corresponds to $t$ by setting $$t' = \lim_{n \to \infty} \phi^{-n}(supp(S^n(t)).$$
In Figure \ref{DPV3} we compare level-$5$ tiles from the DPV (left) and the self-similar tiling (right).
\begin{figure}[ht]
\hspace{1cm}\includegraphics[width=2in]{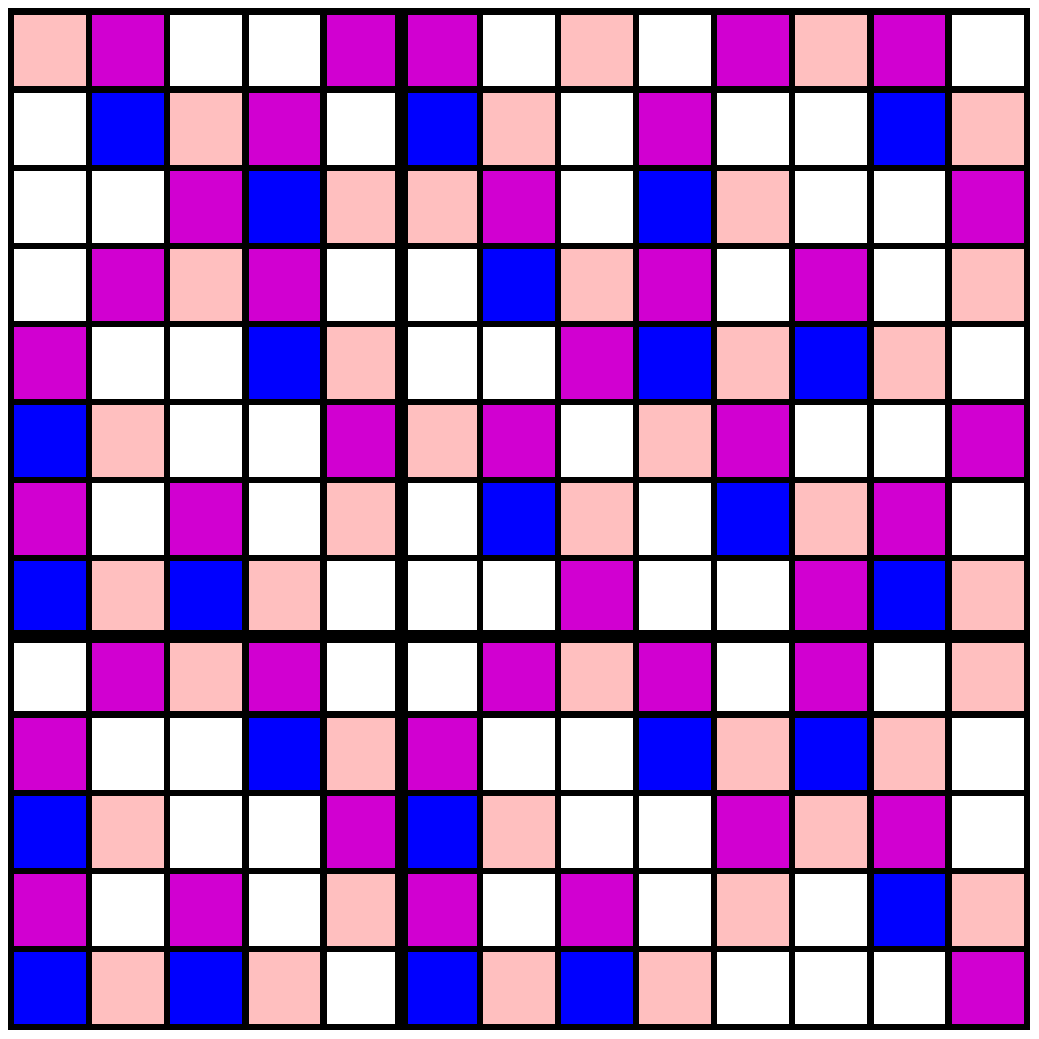} \hspace{1.5cm} \includegraphics[width=2in]{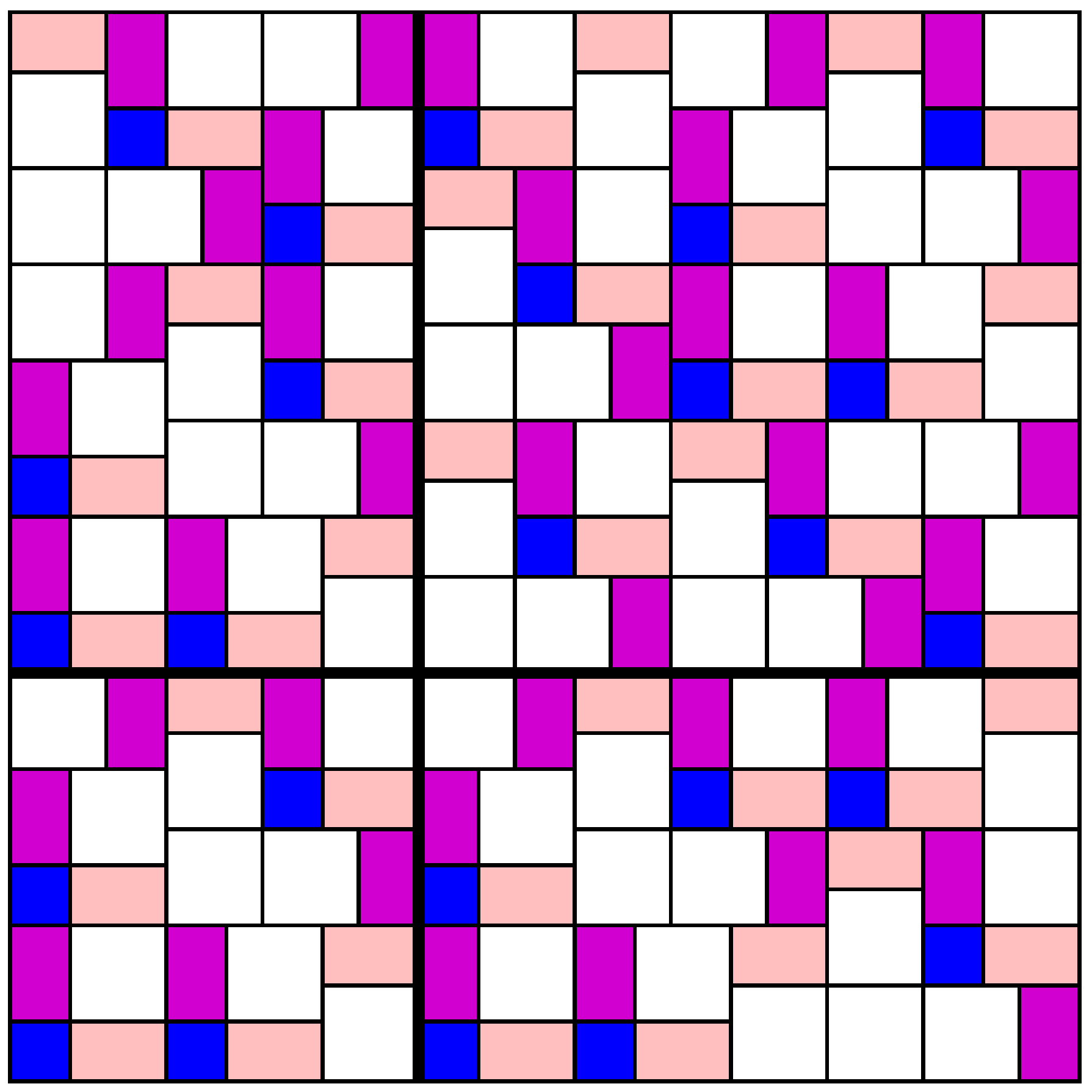}\hspace{1cm}
\caption{Comparing the DPV with the SST of Example \ref{fib.sst}.}
\label{DPV3} 
\end{figure}
\end{ex}

\begin{ex}\label{Rauzy.sst}
The self-similar tiling associated with the Rauzy two-dimensional substitution of Example \ref{Rauzy.cst} has as its
volume expansion the largest root of the polynomial $x^3 - x^2 - 1$.  The three tile
types obtained by the replace-and-rescale method are shown in Figure \ref{Rauzy4}, compared with a large iteration of the substituton.   

\begin{figure}[ht]
\includegraphics[width=5.4in]{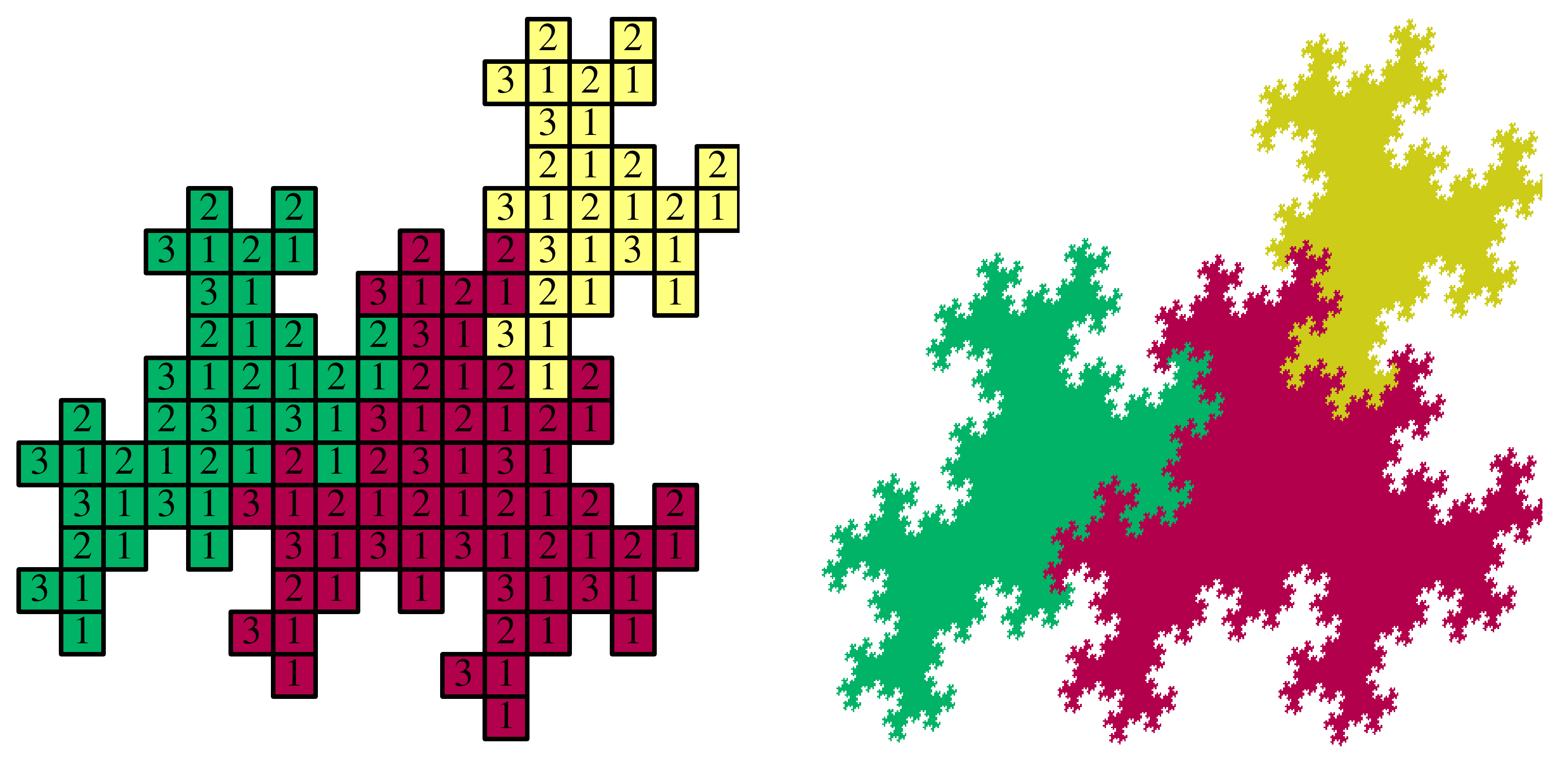}
\caption{A comparison of an iterate with the limiting self-similar tiles.}
\label{Rauzy4} 
\end{figure}
\end{ex}

\subsection{Curious examples}
\label{curiousconnections}
The replace-and-rescale method can produce intriguing results, especially if the substitution is not constructive or not primitive.
We look at the former case in Examples \ref{np.sst} and \ref{trinpsst} and discover that the associated geometric substitution tilings may lose local finiteness.   In Example \ref{Chacon} we consider the latter case to see how a lack of primitivity can impact the geometric substitution; in this case an attempt to ``fix" the situation yields new tiling substitutions that fail to have the expected relationship to one another.

\begin{ex}\label{np.sst}
Applying the replace-and-rescale method to the substitution in Example \ref{np.dpv} produces the inflate-and-subdivide rule of Figure \ref{Non-Pisot2}.
\begin{figure}[ht]
\includegraphics[height=.55in]{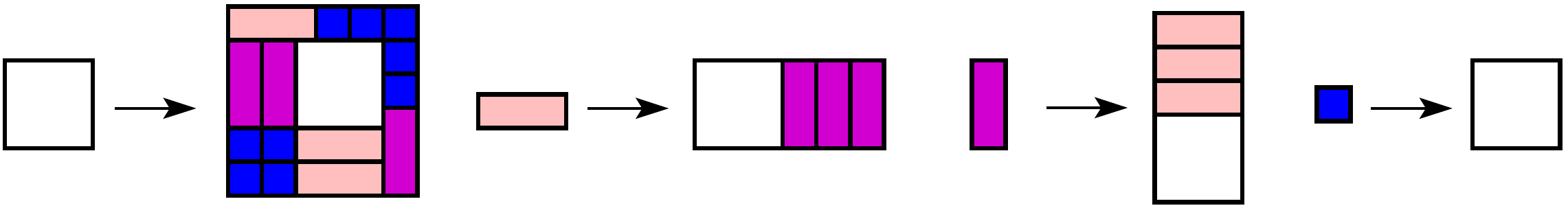}
\caption{The inflate-and-subdivide rule associated with Example \ref{np.dpv}.}
\label{Non-Pisot2} 
\end{figure}
It is proved in \cite{Merobbie} that any tiling admitted by this inflate-and-subdivide rule
does not have finite local complexity since there are two tiles that meet in infinitely many different ways.
(Examination of the tiling on the right of Figure \ref{Non-Pisot3} may convince the reader that this is plausible).  This lack of local finiteness means that the dynamical results found in \cite{Sol.self.similar}, most notably that the system should be weakly mixing, cannot be directly applied.  

\begin{figure}[ht]
\includegraphics[height=3in]{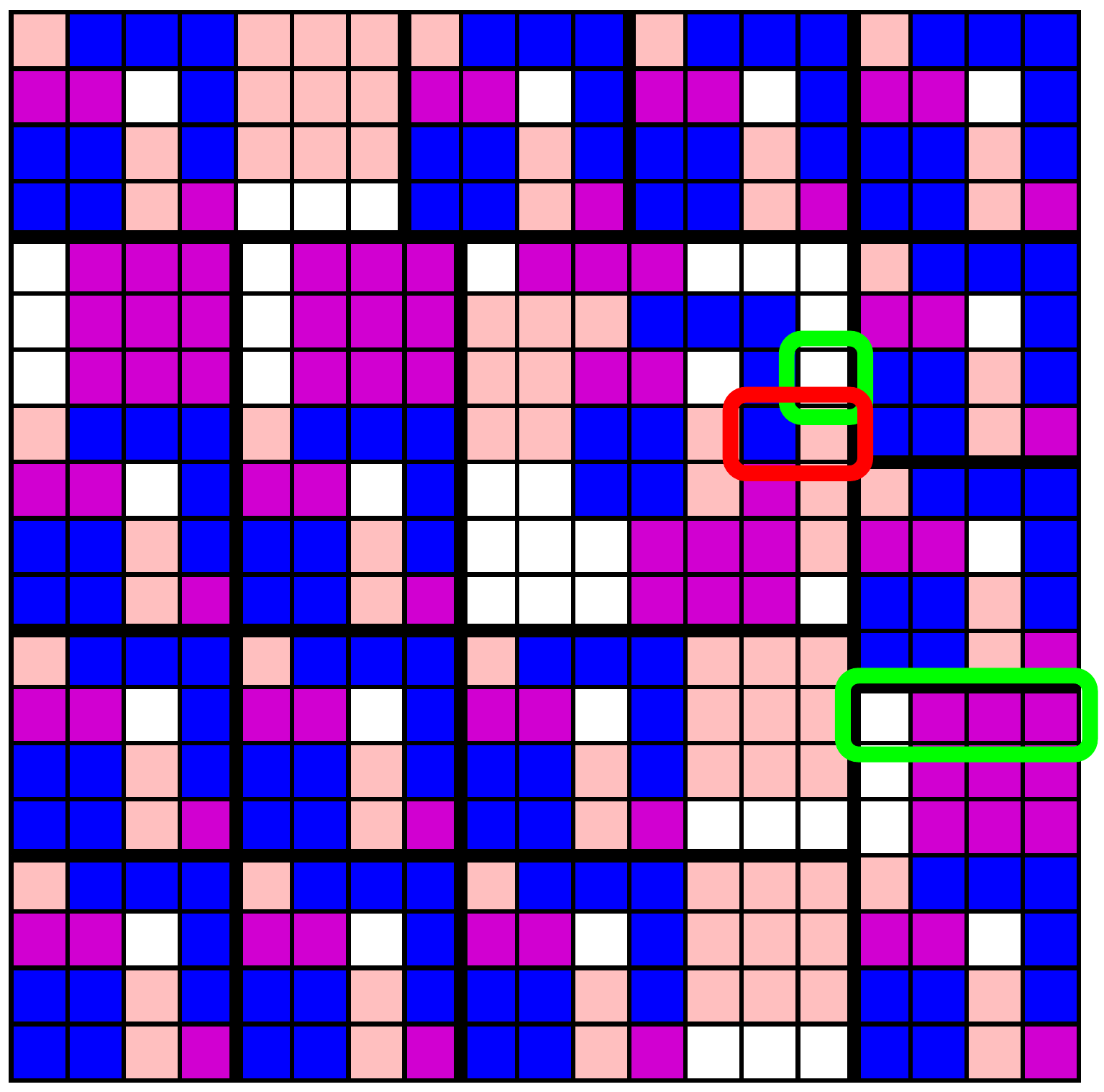} \hfill \includegraphics[height=3in]{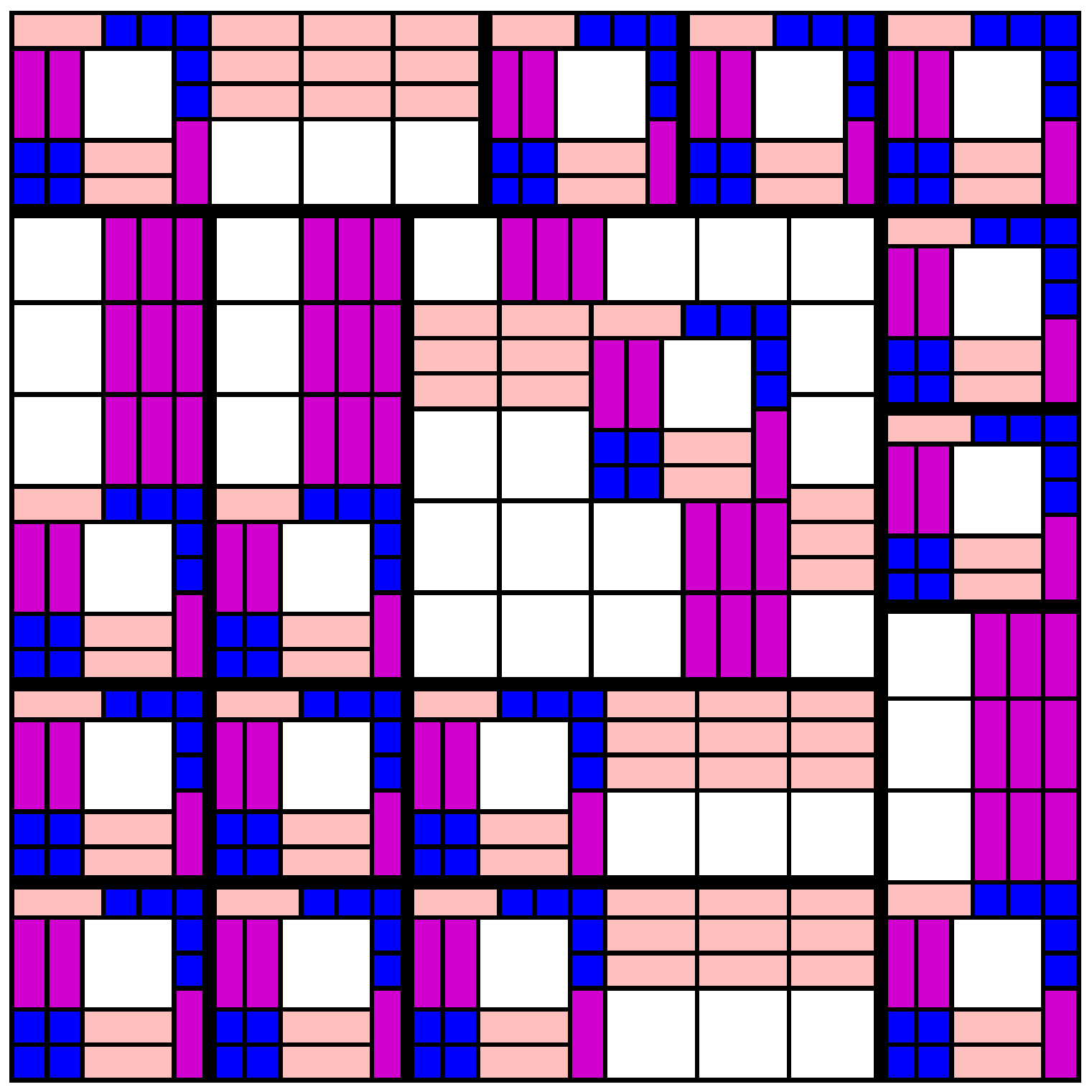}
\caption{Three iterates of the tile of type $1$, with a loss of FLC on the right.}
\label{Non-Pisot3} 
\end{figure}

The loss of finite local complexity happens along arbitrarily long line segments composed of tile edges (proved for this example in \cite{Merobbie} but was found as a necessary condition in general in \cite{Kenyon.rigid}).   As you travel along such a segment, a discrepancy in the number of short tile edges from one side of the line to the other appears; on longer segments this discrepancy increases as more and more short edges pile up on one side than the other.    Because the tile edge lengths are not rationally related, this means that we must keep seeing new adjacencies as the discrepancy  grows.
In the limit one will see infinite {\em fault lines} along which tiles may slide across one another with arbitrary offsets.

The growth of these discrepancies is made possible by the fact that the expansion constant's algebraic conjugate is greater than one in modulus (i.e. it is strongly non-Pisot).   This is also responsible for the fact that original DPV has adjacencies that are ripped apart when substituted, as shown in the left of Figure \ref{Non-Pisot3}.

\end{ex}

\begin{ex}\label{trinpsst}
As in the previous example, the substitution of Example \ref{nptriangles} gives rise to a self-similar tiling that does not have finite local complexity.    In the previous example, fault lines could occur both horizontally and vertically.   In this example, fault lines can occur horizontally, vertically, and diagonally, as one can see from the right side of Figure \ref{Non-pisottriangles2}.   The author has not seen examples allowing fault lines in more than three distinct directions.
\begin{figure}[ht]
\includegraphics[width=3in]{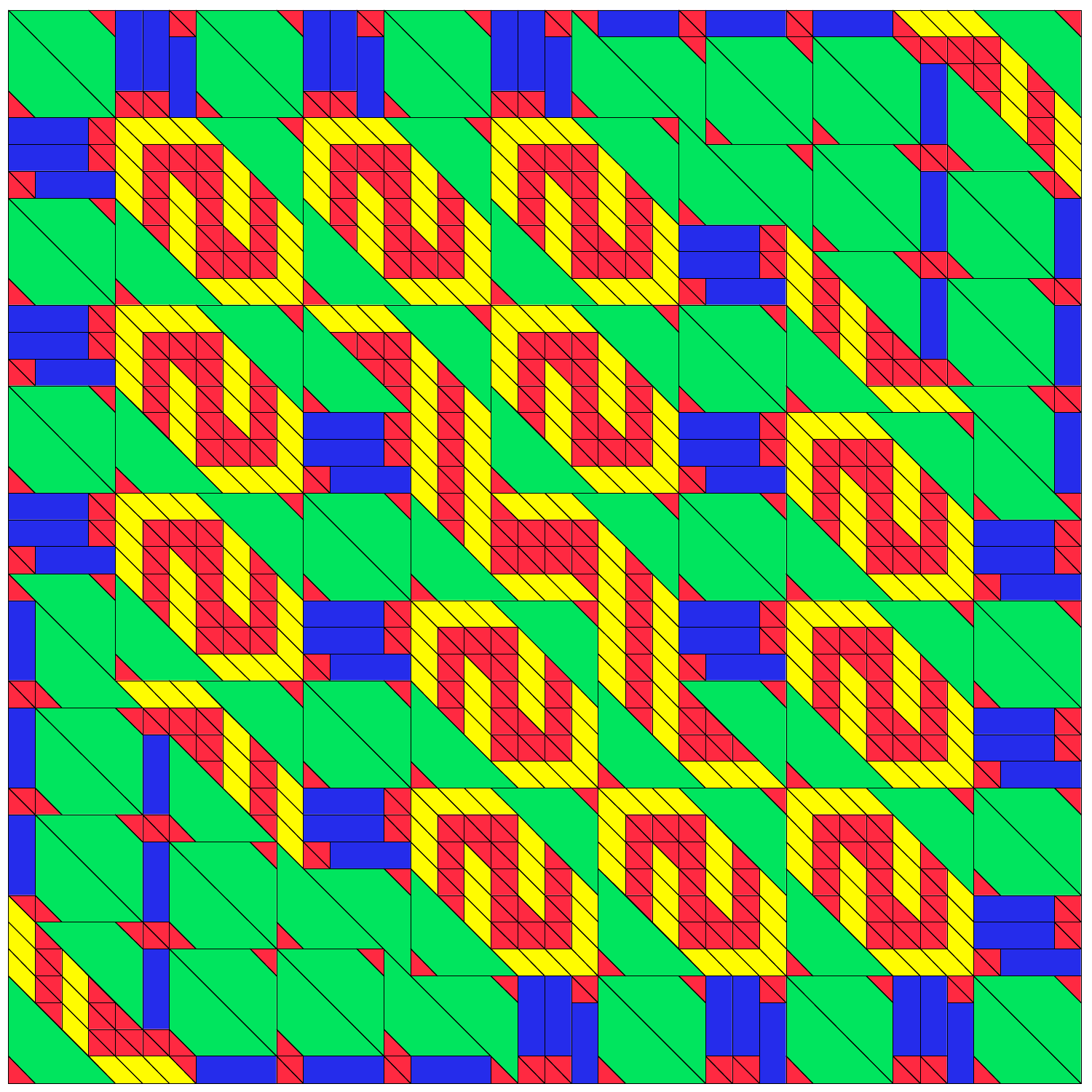} \hfill \includegraphics[width=3in]{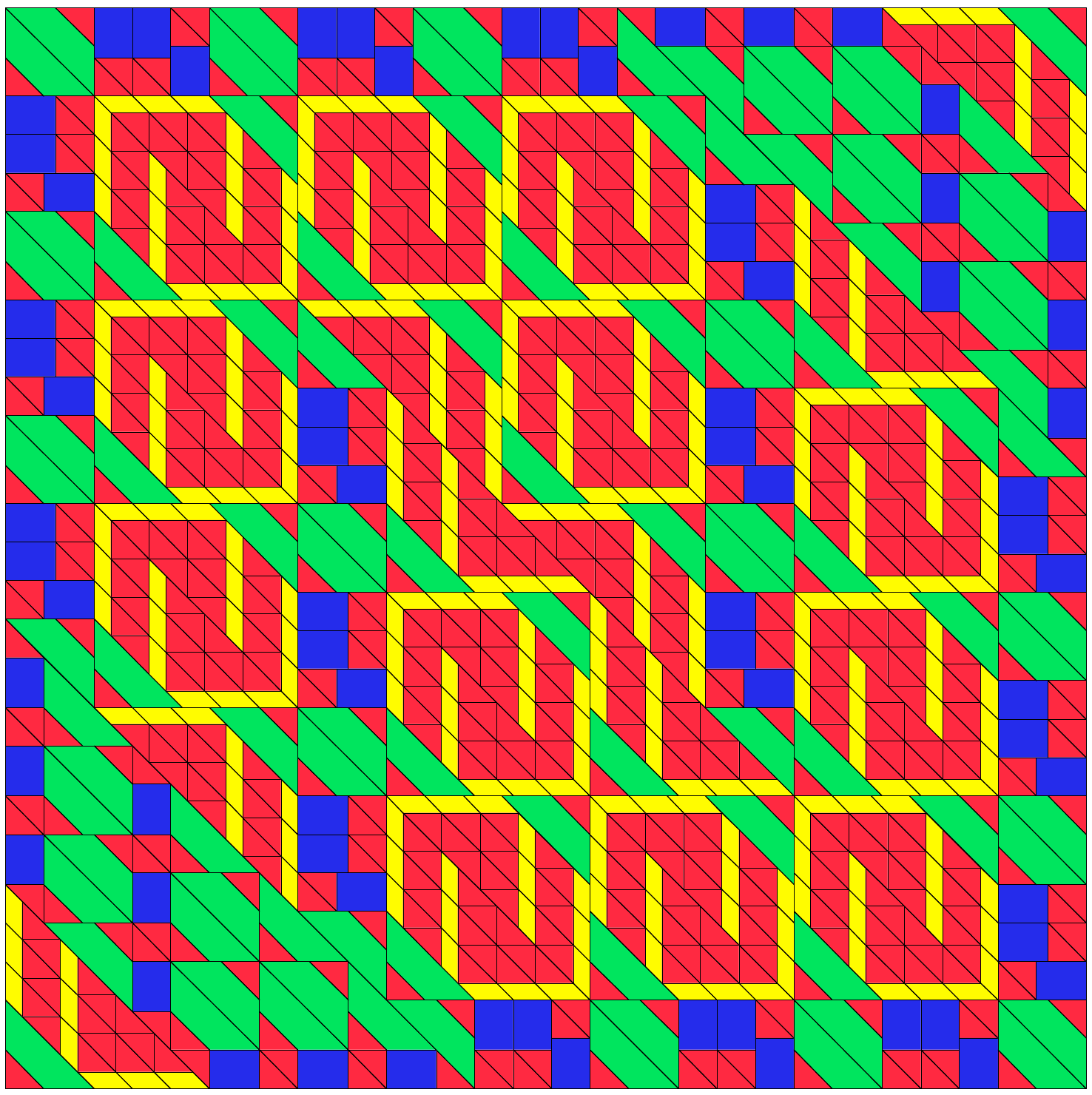}
\caption{Four iterations of a square formed by two red triangles.  Fault lines occur in three directions in the tiling on the right.}
\label{Non-pisottriangles2} 
\end{figure}

\end{ex}

\begin{ex}\label{Chacon}
A famous one-dimensional dynamical system is given by the Chacon cut-and-stack construction, which provided the first example of a weakly but not strongly mixing system (see \cite{Fogg}, p. 133 for a synopsis of the results in the one-dimensional case).
The cut-and-stack system can be recoded by the symbolic substitution $a \to a a b a, b \to b$, and Figure \ref{Chaconsq} shows a DPV substitution based on this.
\begin{figure}[ht]
\includegraphics[width=4.25in]{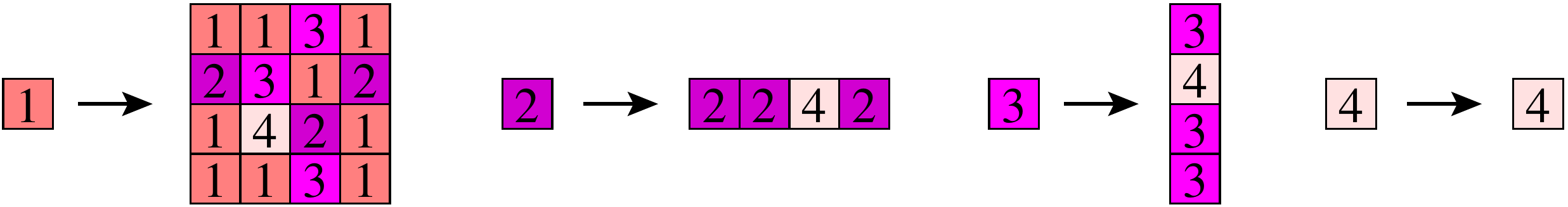}
\caption{A Chacon DPV substitution.}
\label{Chaconsq} 
\end{figure}

Another tiling version of
the construction, shown in Figure \ref{Chacon1}, is analyzed from a dynamical systems viewpoint in \cite{Rob.Park} and put in the context of combinatorially substitutive tilings in \cite{Mycstpaper}.   The four prototiles used in those works are not square, but are  a rescaling of the supports of the level-1 tiles in Figure \ref{Chaconsq}.   
\begin{figure}[ht]
\includegraphics[width=4.25in]{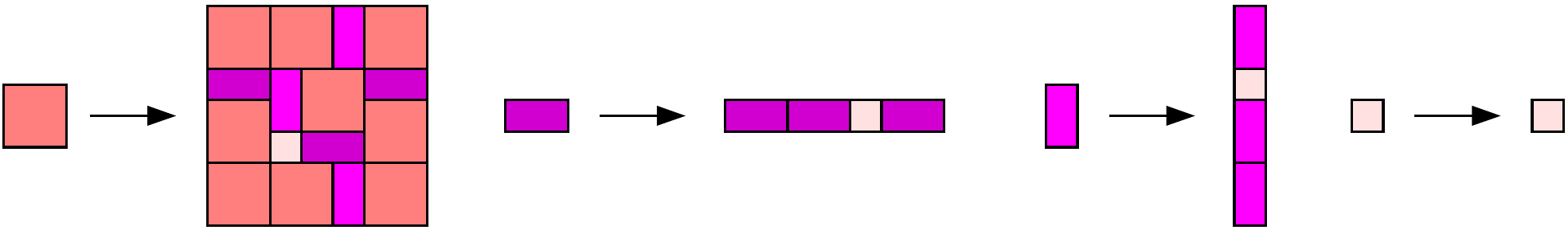}
\caption{The Chacon cut-and-stack construction found in \cite{Rob.Park}.}
\label{Chacon1} 
\end{figure}

This substitution is not primitive since no matter how many times
we substitute the three small tiles, they will never contain the large square.  Because of this we cannot obtain a meaningful self-similar tiling directly using the replace-and-rescale method: the replacements of all but the large square will have volumes that go to zero under rescaling, leaving us with only the first tile and a trivial substitution.

There is a way to recode the system into a primitive one, producing the prototiles shown to the right of the arrows in Figure \ref{Chacon3}.   By referring to Figure \ref{Chacon2}, the reader may be convinced that knowing the surroundings of a particular tile is enough to decide unambiguously with which new prototile to replace it. 
\begin{figure}[ht]
\includegraphics[width=5in]{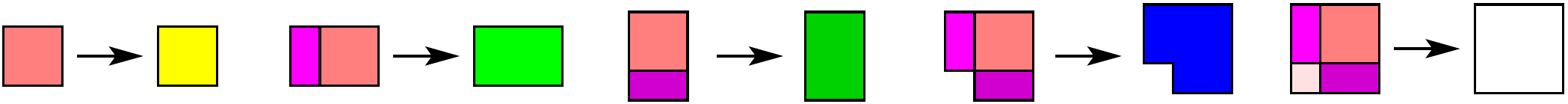}
\caption{A new tile set that makes a primitive Chacon substitution.}
\label{Chacon3} 
\end{figure}
\begin{figure}[ht]
\includegraphics[width=2in]{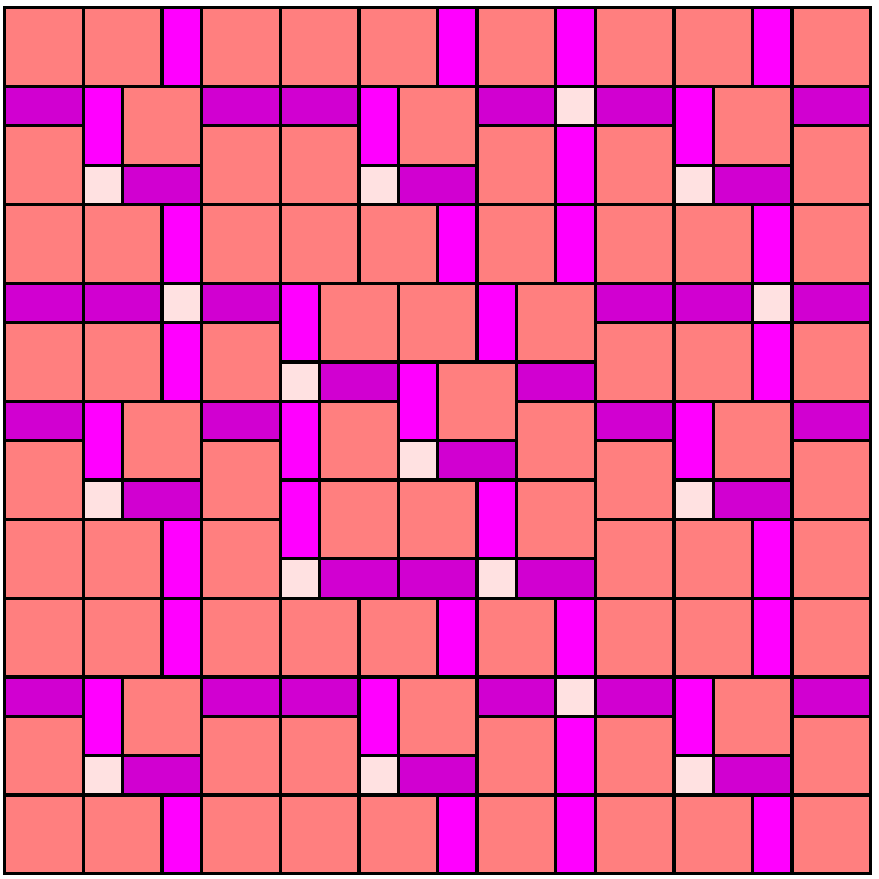} \hspace{1cm}\includegraphics[width=2in]{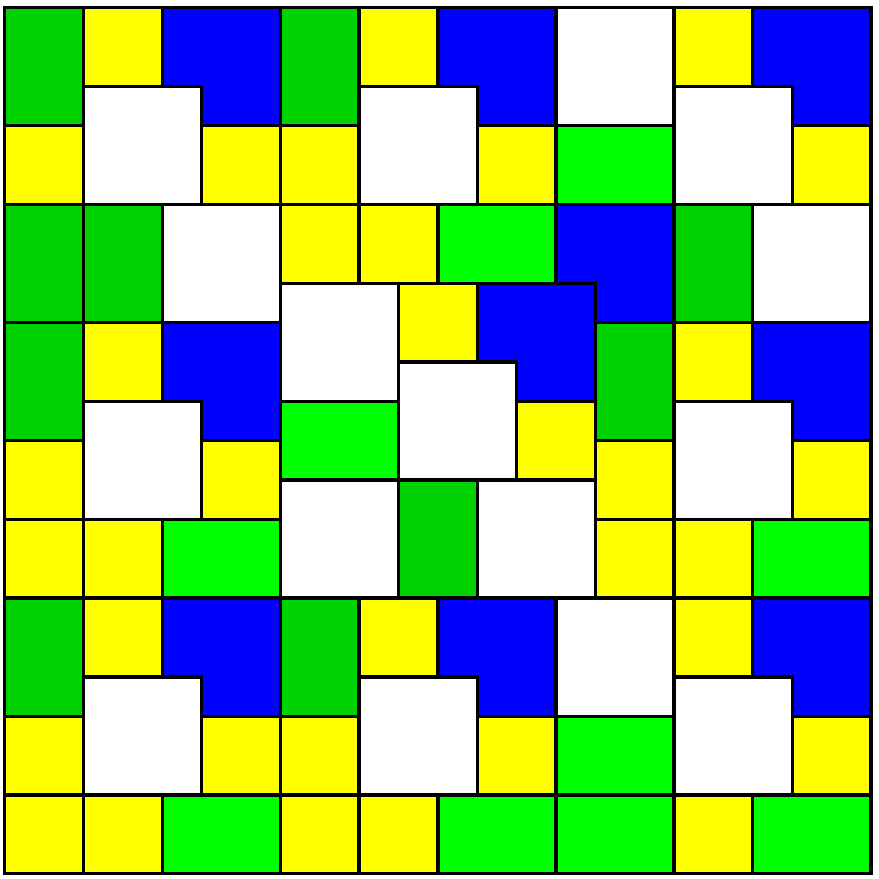}
\caption{Chacon level-two tiles, nonprimitive and primitive versions.}
\label{Chacon2} 
\end{figure}
The nonprimitive Chacon substitution turns into the primitive one of Figure \ref{Chacon4}.    
\begin{figure}[ht]
\includegraphics[width=3.25in]{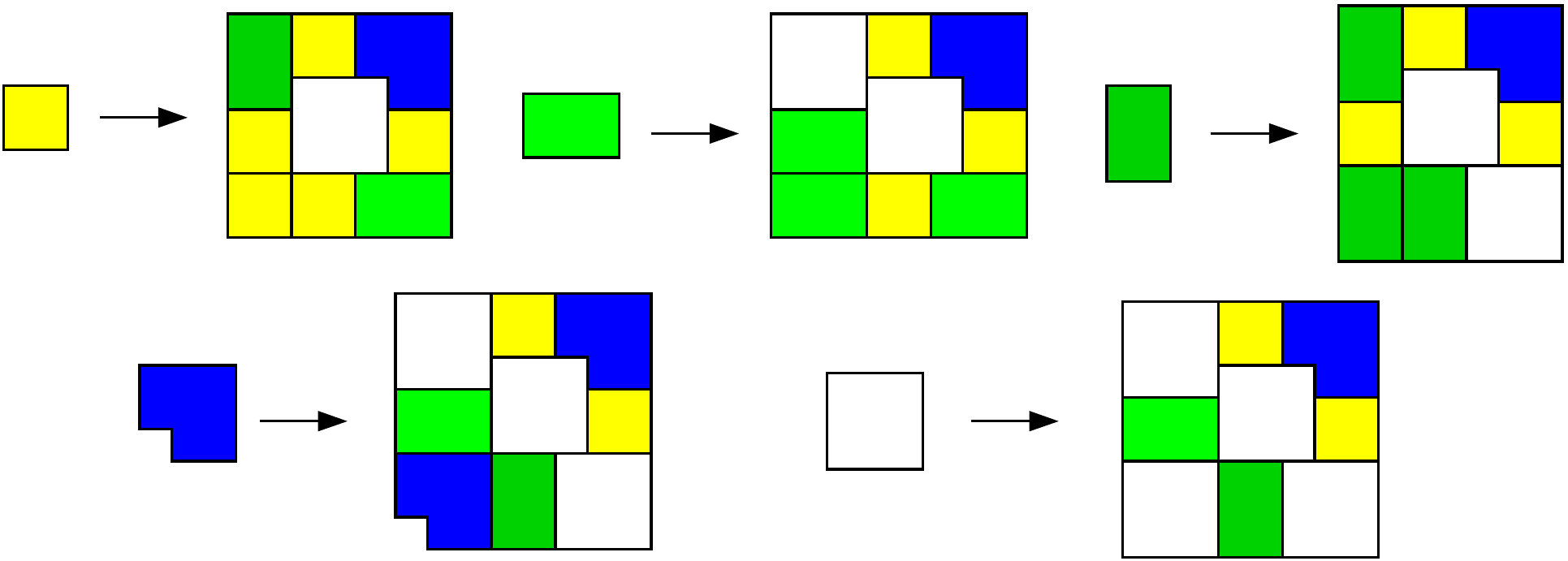}
\caption{The Chacon primitive substitution.}
\label{Chacon4} 
\end{figure}

Because there is a locally defined map taking tilings admitted by the nonprimitive substitution to tilings admitted by the primitive substitution, and vice versa, the tiling spaces are considered {\em mutually locally derivable}.    This means that the dynamical systems are equivalent in the sense of ``topological conjugacy" \cite{Mypaper1}, thus important dynamical features are preserved.
One such dynamical feature, proved in \cite{Rob.Park}, is that the dynamical system under the action of $\Z^2$ is weakly mixing.

The larger eigenvalue of the substitution matrix of the Chacon primitive substitution is 9, and it is not difficult to see that the length expansion is governed by powers of 3.  The replace-and-rescale method produces a prototile set of five congruent squares; the inflate-and-subdivide rule is shown in Figure \ref{Chacon5}.
\begin{figure}[ht]
\includegraphics[width=3.1in]{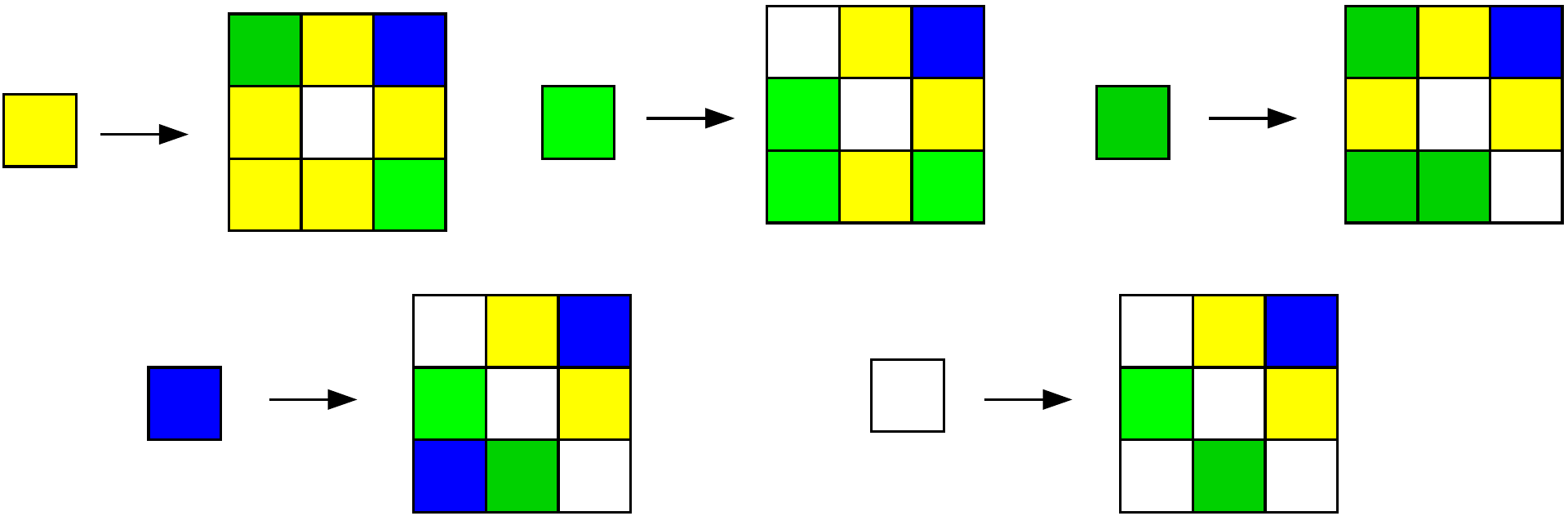}
\caption{The inflate-and-subdivide rule associated to the Chacon primitive substitution.}
\label{Chacon5} 
\end{figure}

What makes this example curious is that the dynamical systems of the combinatorially substitutive tiling and its associated self-similar tiling are distinctly different.
One can use results from \cite{Sol.self.similar} to show that under the $\R^2$ action, the self-similar tiling dynamical system is not weakly mixing.   An embedded $\Z^2$ action would therefore also fail to be weakly mixing. 
This stands in contrast to the weakly mixing $\Z^2$ action proved in \cite{Rob.Park} when the substitution is only combinatorial.
One can see that the systems are ``misaligned" by considering Figure \ref{Chacon6} and comparing the location of the red circle in each substitution, which for $n = 0, 1,$ and $2$ represents the central level-$n$ tile within its level-$(n+1)$ tile.  One can check that as
$n$ grows without bound so does the distance between the red dots.   That is, if any corners of the level-$n$ tiles are lined up, the red dots will move further and further away from one another!   

\begin{figure}[ht]
\includegraphics[width=6.0in]{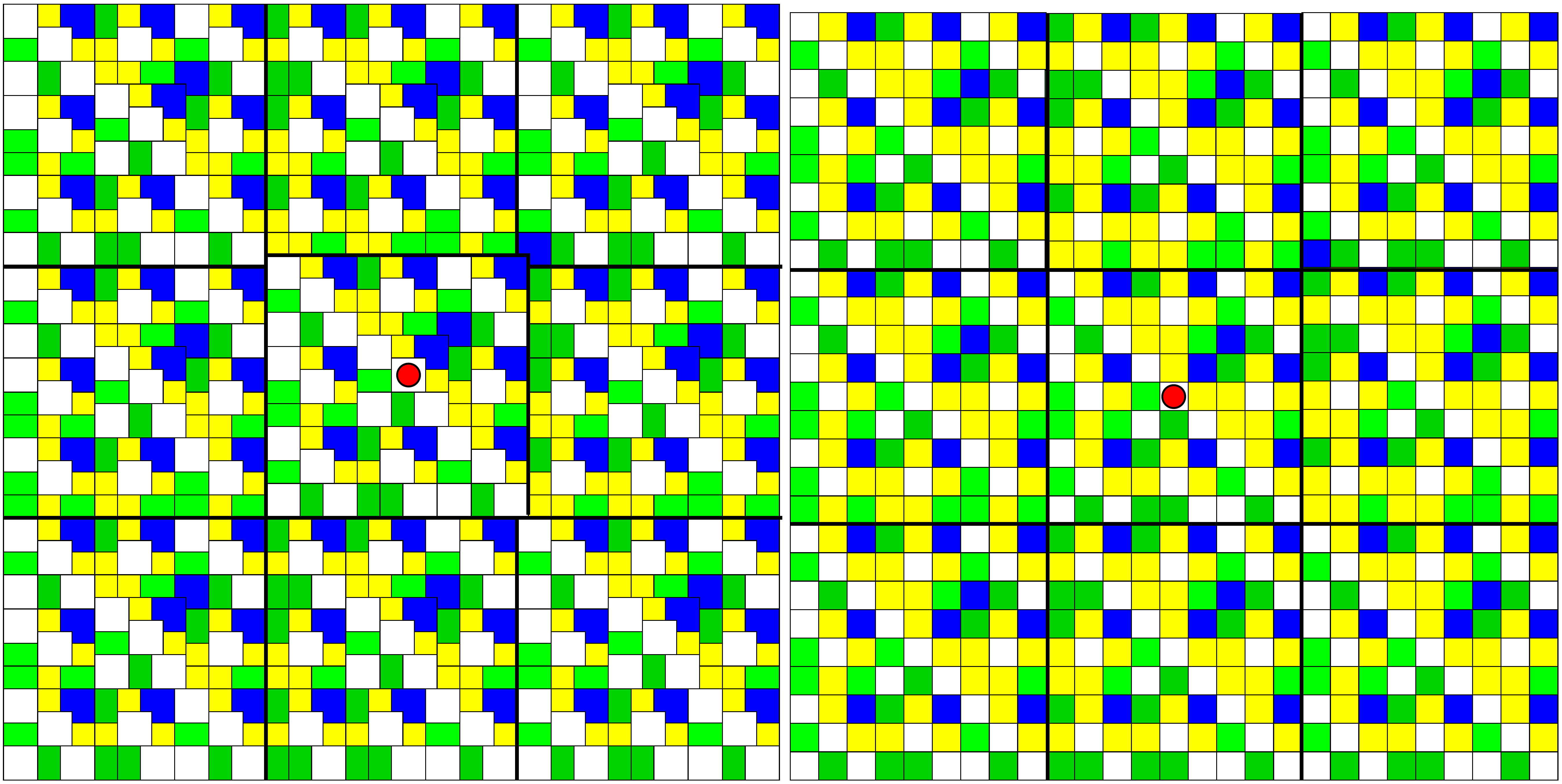}
\caption{Comparing three iterations of the white square.}
\label{Chacon6} 
\end{figure}

\end{ex}

\section{Some areas of research}
\label{Results.questions}
Substitution tilings are being studied from topological, dynamical, physical, combinatorial, and other perspectives, often in conjunction with one another.   In this section we will briefly outline areas of current interest and possible questions for future study.

\subsection{Geometric tiling substitutions}
\label{Geom.results.questions}
Tilings make good models for the atomic structure of crystals and quasicrystals, and perhaps the most exciting work on them is being done at the intersection of physics and topology.  Methods for investigating certain tiling spaces via $C^*$-algebras have been developed
and are nicely summarized in \cite{Kellendonk.Putnam}.   The types of tilings that are most easily evaluated this way are self-similar tilings and tilings generated by the projection method.    (Some tilings, such as the Penrose and the ``octagonal" tilings, fall into both categories).
The K-theory of these $C^*$-algebras are of interest to both mathematicians and physicists.  The possible energy levels of electrons in a material modeled by a tiling determine gaps in the spectrum of the associated Schr\"odinger operator.   The K-theory gives a natural labeling of the spectral gaps, thus providing theoretically relevant physical information (see \cite{Bellissard.etal} for a detailed discussion of this branch of study).   It is believed that there may be additional physical interpretations for K-theory and other topological invariants of tiling spaces.


There is more promising topological work being done as well.    For instance, it has been shown in varying degrees of generality (see \cite{Sadun_limit} and references therein) that FLC tiling spaces are inverse limits.   Successful efforts to compute the homology and cohomology of tiling spaces, and to connect these results to K-theory, have been plentiful.    A nice summary of the current state of the art, along with the discovery of torsion and its ramifications, appears in \cite{Gahler.Hunton.Kellendonk} with a primary emphasis on ``canonical" projection tilings.  An informal discussion of the connections between some physical and mathematical problems appears in \cite{Sadun.top.survey}, with a focus on recent progress in the cohomology of tiling spaces.     Included is a summary of the work in \cite{Clark-Sadun} involving cohomological analysis of the deformations of tiling spaces.  An important question is the extent to which the homology and cohomology of tiling spaces has physical interpretations.  

Almost all of the existing literature on the topology of tiling spaces makes the assumption of local finiteness.   This is, after all, an appropriate restriction, given that the model of atomic structure requires tiles (atoms) to fit together in only a finite number of ways.   
However, examples exist of geometric tiling substitutions that result in non-FLC tilings, for example Danzer's triangular tilings \cite{Danzer}
and the tiling from \cite{Merobbie}, which is easily generalized using the methods of Section \ref{connections.section}.    In \cite{Kenyon.rigid}, Kenyon was the first to consider the conditions under which a tiling of $\R^2$ with finitely many tile types can lose local finiteness: the tile boundaries must contain circles or arbitrarily long line segments, thus substitution tilings without finite local complexity have fault lines along which tiles can slide past one another.  In \cite{Melorenzo}, the cohomology of a highly restricted class of non-FLC substitution tilings
was successfully computed, and it was shown that each fault line leaves a sort of signature in the cohomology in dimension 3, even though the tilings are two-dimensional.
It is a topic of current interest to understand the topology of tiling spaces without finite local complexity.

At the intersection of mathematical physics and dynamical systems is the connection between the diffraction spectra of quasicrystalline solids and the dynamical spectra of the tilings that model them.   The fact that these spectra are related at all is first established in \cite{Dworkin}, while the mathematical description of the diffraction spectral measure is given sound theoretical footing in \cite{Hof}.    
Much of the work to date has centered around discrete point sets called Delone sets, which can be thought of as locations of molecules and which can be converted into tilings in a few different ways.     
Ever since Schectman et. al. \cite{Schechtman} discovered quasicrystals in a laboratory experiment, people have been trying to figure out which
Delone sets are ``diffractive" in that their spectra exhibit sharp bright spots.   Mathematically it is interesting to ask when the spectra consists only of sharp bright spots, i.e. when it has ``pure point spectrum".  More precisely, one defines a
{\em spectral measure} which can be broken into pure point, singular continuous, and absolutely continuous pieces with respect to Lebesgue measure.
Great progress has been made for ``model sets" (obtained by a generalized projection method), and for Delone sets generated by substitutions; a current synopsis of the state of the art appears in \cite{Moody.survey}.
 It is now known that for certain locally finite Delone sets the notions of pure point dynamical and pure point diffraction spectra coincide \cite{Lee.Moody.Solomyak}.  This was generalized in \cite{Baake.Lenz} in a measure-theoretic setting which allows for a lack of local finiteness.  
The question of whether certain substitution systems consist of model sets can be investigated by looking for ``modular coincidences"; \cite{Frettloh.Sing} has an algorithm and many examples, which build upon the work in \cite{Lee.Moody.Solomyak}.
Questions remain regarding the connections when there is any continuous portion of the spectral measures.
The dynamical spectra of specific geometric tiling substitutions have been studied (\cite{Sol.self.similar,Mytoppaper} and others) but are not completely understood. 

Related to the study of tilings and model sets is a question in dynamical systems theory.  For one-dimensional symbolic substitutions, it is sometimes possible to find a ``geometric realization" of the substitution.   A formal definition appears in \cite{Queffelec}, p. 140, but the idea is that a geometric realization is a geometric dynamical system, (such as an irrational rotation of the circle), which encodes the system via partition elements.
For example, the Fibonacci substitution sequences can be seen to code, in an almost one-to-one fashion, addition by  $1/\gamma$ on a one-dimensional torus, where $\gamma$ is the golden mean (see \cite{Fogg}, p. 199 for details).   Orbits in this geometric realization ``look like" one-dimensional tilings.   Several more examples are given in \cite{Fogg}, p. 231.    If a tiling dynamical system arises from a model set, then it can be seen as a geometric realization.  
For example, it is shown in  \cite{Rob.penrose}  that the Penrose rhombic tiling dynamical system is an almost one-to-one extension of an irrational rotation on a 4-dimensional torus.   In general, we do not know when substitution tilings have geometric realizations.

\subsection{Combinatorial tiling substitutions}
\label{Comb.results.questions}
This paper provides the full extent, to the author's knowledge, of known classes of combinatorial tiling substitutions.   All of the examples in Section \ref{comb.section} are obtained by various means from one-dimensional symbolic substitutions.  What other mechanisms exist for generating combinatorial substitutions?  Is there a method for obtaining non-geometric substitutions from geometric ones?    It seems clear that there should be a multitude of other examples waiting to be discovered, and finding them is of paramount importance.

Combinatorial tiling substitutions have hardly begun to be studied from the dynamical systems viewpoint.   In analogy with the self-similar case and many of its generalizations, we would like to investigate basic ergodic-theoretic properties such as repetitivity, unique ergodicity, and recognizability.   This program was carried out in \cite{Hansen} on a restricted class of two-dimensional symbolic substitutions of non-constant length for which ``standard" techniques could be applied.
Unfortunately, these techniques  do not necessarily work in the non-geometric case.   The crucial missing piece is that 
the substitution rule cannot always be seen as an action from the tiling space to itself:  the substitution can be applied only to level-$n$ tiles, not to entire tilings.   
Many combinatorial tiling substitutions do not extend to maps of the tiling space in a canonical way, and it is unclear whether (or when) any of them do.  New methods will need to be devised to tackle even the most basic questions in the dynamical systems and ergodic theory of combinatorial substitution tilings.

 A closely related concept, essential to many standard arguments, is whether the substitution map can be locally ``undone" so that one can detect the level-1 tile in a given region without requiring infinite information about the tiling.    This is called {\em recognizablility} in the sequence case and the {\em unique composition property} in the self-similar tilings case.  When the substitution acts as a continuous map on the tiling space, unique composition is equivalent to the substitution map being invertible.   In the event of non-periodicity, recognizability and unique composition were proved in \cite{Mosse} and \cite{Sol.u.comp}, respectively.   Although the substitution map may not make sense on tiling spaces in the non-geometric case, the notion of unique composition still does, and a natural conjecture is that combinatorial substitutions possess it whenever they are non-periodic.

In one dimension, there is great interest in the theory of ``combinatorics on words" (see Part I of \cite{Fogg} for an extensive exposition).   In this theory, one considers finite blocks of letters and investigates how often they appear, and in what combinations, within sequences.  Substitution sequences are particularly fertile for this type of study.  The {\em complexity} of a sequence is a function $p(n)$ telling how many words of length $n$ exist in the sequence; this can be used to compute the topological entropy \cite{Fogg}, p. 4.   Any non-periodic sequence with minimal complexity is called {\em Sturmian}, the classic example being the sequence given by the Fibonacci substitution.  One can read about the numerous consequences of being Sturmian in Chapter 6 of \cite{Fogg}.    The notion of complexity can be generalized to higher dimensions and some results exist in this direction (see \cite{Berthe.Tijdeman} and references therein).  Combinatorial substitutions such as the Rauzy substitution of Example \ref{Rauzy.cst} are a natural place to look for examples of low-complexity sequences.

Some problems that have been at least partially resolved for geometric substitutions are still open for combinatorial ones.   For instance it is completely unclear whether there should exist matching rules which force tiles to fit together as prescribed by combinatorial substitutions.   Would it be possible to use the matching rules for their associated self-similar tilings, which we know exist by \cite{Goodman-Strauss}, to find them?   Another question is, since the connection to the atomic structure of solids is an important motivation for the study of tiling spaces, can we identify the diffraction or dynamical spectrum of combinatorial substitution?    It is known that the dynamical spectrum of the Chacon substitution is trivial since it is weakly mixing \cite{Rob.Park}, and following \cite{Sol.self.similar} it is reasonable to conjecture that DPV substitutions without Pisot expansions might also be weakly mixing.   For combinatorial substitutions, is the spectrum largely dependent on the Perron eigenvalue of the substitution matrix, as it is in the geometric case \cite{Sol.self.similar}?  Or is the situation like the one-dimensional symbolic case, which is also highly sensitive to the combinatorics of the substitutions?

\subsection{Connections between the constant and non-constant length cases}
\label{Conn.results.questions}
The first open question is, when does a combinatoral tiling substitution give rise to a reasonable geometric one?
We have already seen that the non-primitive Chacon substitution of Example \ref{Chacon} does not.  There must be substitutions for which the limit in the replace-and-rescale method does not exist, or produces topologically unpleasant tiles.    In fact, it is unclear exactly how the replace-and-rescale method ought to properly be applied:  determining the appropriate linear expansion map is problematic for at least two reasons.    First,  the combinatorial substitution might be encoding an inflate-and-subdivide rule that does not inflate as a similarity.   This means that knowing the volume expansion would not tell us the appropriate length expansions.  Second, if the linear map is a similarity, there may be some rotation inherent in the combinatorial substitution that would need to be expressed in the linear map, as in the Rauzy substitution of Example \ref{Rauzy.cst}.   In the best circumstances we could hope to find conditions under which the expansion can be found, the limiting tiles are topologically ``nice", and a proper inflate-and-subdivide rule exists.

It is interesting to consider the relationship between the dual graphs of a  combinatorial substitution and its associated geometric substitution (if it exists).   It is clear from Figure \ref{DPV3} that the dual graphs must always have the same labeled vertices, but the edge and facet sets do not seem to bear a consistent relationship to one another.  In the case of Example \ref{fib.sst} the edge set of the DPV's dual graph is contained in that of its self-similar tiling, but this is not true in general.
Since the unlabeled dual graphs are not isomorphic, there is no homeomorphism of the plane taking one to another (see \cite{G-S}, p. 169).   Can an understanding of the combinatorial properties of one tiling still give us insight into the other?

In one-dimensional symbolic dynamics, the Curtis-Lyndon-Hedlund Theorem (see \cite{Lind.Marcus}) states that homeomorphisms between symbolic dynamical systems are equivalent to local maps called  ``sliding block codes".   A sliding block code transforms one sequence into another  element by element, deciding what to put in the new sequence by looking in a finite window in the old one.    Similarly, one tiling can be transformed locally into another; if the process is invertible the tilings are mutually locally derivable.
It has come to light that there is no Curtis-Lyndon-Hedlund theorem for tiling dynamical systems \cite{Petersen.oned, rasad}.   
Using the basic method of \cite{rasad}, one can show that Example \ref{fib.dpv} and the associated self-similar tiling of Example \ref{fib.sst} have topologically conjugate dynamical systems without the possibility of mutual local derivability.   We conjecture that in the Pisot case, combinatorial substitutions have topologically conjugate dynamical systems with their geometric counterparts.  In general one would not expect the conjugacy to be through mutual local derivability.

Our final question takes note of the fact that the dynamical relationship between substitution sequences and self-similar tilings of the line is especially subtle.  
On a sequence space there is a $\Z$-action; passing to a tiling by choosing tile lengths provides a natural action by $\R$ called a ``suspension".   
Surprisingly, the continuous action of the tiling space is probably better understood than the discrete action!  For instance, the presence or absence, and nature of, eigenvalues of the tiling dynamical system can be understood in terms of the expansion constant along with certain geometric information \cite{Sol.self.similar}.   This situation is far more complicated in the symbolic case and the interested reader should see \cite{Fogg}, Section 7.3 for a synopsis.
Also, topological conjugacies between different suspensions have been thoroughly considered in \cite{Clark-Sadun1D}, where it is seen that the eigenvalues of the substitution matrix play a critical role.
We can consider tilings such as those in Figure \ref{Non-Pisot3} as being suspensions of the same sequence in $\Z^2$, and ask similar questions about their spectra and topological properties.   More generally, we can consider tilings such as those in Figure \ref{Non-pisottriangles2} as being suspensions of the same labeled graph.
This perspective yields an interesting set of problems at the intersection of dynamics and combinatorics.

\bibliographystyle{amsplain}

 \end{document}